\documentclass{amsart}
\usepackage{bbm}
\usepackage{dsfont}
\usepackage{relsize}
\usepackage{amsthm}
\usepackage[breaklinks]{hyperref}
\usepackage{blindtext}
\usepackage{amsmath,amssymb}
\usepackage{autonum}
\usepackage[usenames,dvipsnames,svgnames]{xcolor}
\usepackage{mathrsfs}
\usepackage{csquotes}
\usepackage{mathtools}
\usepackage[]{hyperref}
\hypersetup{colorlinks=true,linkcolor=WildStrawberry,citecolor=cyan}
\newtheorem{thm}{Theorem}[section]
\newtheorem{cor}[thm]{Corollary}
\newtheorem{lem}[thm]{Lemma}

\theoremstyle{hypothesis}

\theoremstyle{definition}
\newtheorem{defn}[thm]{Definition}
\theoremstyle{remark}
\newtheorem{rem}[thm]{Remark}
\theoremstyle{example}

\numberwithin{equation}{section}

\renewcommand{\a}{\alpha}
\renewcommand{\b}{\beta}
\renewcommand{\c}{\gamma}

\makeatletter
\renewcommand\paragraph{\@startsection{paragraph}{4}{\z@}%
	{-2.5ex\@plus -1ex \@minus -.25ex}%
	{1.25ex \@plus .25ex}%
	{\normalfont\normalsize\centering\bfseries}}
\makeatother
\begin{document}
	\title[Mixed fully nonlinear local and nonlocal sub-elliptic operators in Heisenberg group]
	{MIXED FULLY NONLINEAR LOCAL AND NONLOCAL SUB-ELLIPTIC OPERATOR IN HEISENBERG GROUP: EXISTENCE \& REGULARITY OF SOLUTIONS }
	\author[P.\,Oza,\,\, J.\,Tyagi]
	{ Priyank Oza,\,\,Jagmohan Tyagi }
	\address{Priyank\,Oza \hfill\break
		Indian Institute of Technology Gandhinagar \newline
		Palaj, Gandhinagar Gujarat, India-382055.}
	\email{priyank.k@iitgn.ac.in, priyank.oza3@gmail.com}
	\address{Jagmohan\,Tyagi \hfill\break
		Indian Institute of Technology Gandhinagar \newline
		Palaj, Gandhinagar Gujarat, India-382055.}
	\email{jtyagi@iitgn.ac.in, jtyagi1@gmail.com}
	\thanks{Submitted \today.  Published-----.}
	\subjclass[2010]{Primary 35A01, 35J60, 35R03, 35D40, 47G20; Secondary 45K05}
	\keywords{Nonlocal and local operators, Partial differential equations on the Heisenberg group, Pucci's extremal operator, Integro-PDE, viscosity solutions, Perron's method}
	\begin{abstract}
		We establish the comparison principle, existence and regularity of  solutions to the following problem concerning the mixed operator:
		\begin{align}\label{eq 0.1}
			\begin{cases}
				\alpha\mathcal{M}^+_{\lambda,\Lambda}\big(D^2_{\mathbbm{H}^N,S}u\big)-\beta(-\Delta_{\mathbbm{H}^N})^su=f &\text{in } \,{\Omega},\\
				u=g &\text{in } \,\mathbbm{H}^N\setminus\Omega,
			\end{cases}	
		\end{align}
		where $\mathcal{M}_{\lambda,\Lambda}^+$ is the extremal Pucci's operator and $(-\Delta_{\mathbbm{H}^N})^s$ denotes the fractional sub-Laplacian on Heisenberg group.
	\end{abstract}
	
	\maketitle
	
\section{Introduction}
Recently, mixed operators have been studied in the Euclidean setting by several researchers. These operators occur naturally, for instance, in the study of plasma physics \cite{Blaze}, population dynamics \cite{Dipierro} and many other places. One may see \cite{Barles2, Biagi, Biagi2, Biagi3, Anup mixed, Chen, Chen 2, Gradient} for the existence and related qualitative questions on mixed operators in the Euclidean framework. In those works, the local term is  second order linear elliptic operator and nonlocal term is the fractional Laplacian. Mou \cite{Mou} studied the regularity theory for Hamilton-Jacobi-Bellman-Isaacs type integro-PDEs. Modasiya and Sen \cite{Modasiyaaa} investigated the boundary regularity results for fully nonlinear mixed local-nonlocal equations. There are several articles, where authors have established the existence of solution to  Dirichlet problem using the classical Perron's method, see \cite{bard,User,Crand,Ishii}.  One may also consult \cite{Lara, Mou1} and the reference therein for more related results.

On the other hand, there has been significant works for  $\Delta_\infty$  (infinity Laplace) equations and degenerate equations in non-Euclidean settings, for e.g., Heisenberg group, and more generally on sub-Riemannian, Carnot-Carath\'eodory spaces. Bieske \cite{Bieske} studied infinite harmonic functions in the Heisenberg group using the notion of viscosity solutions. Wang \cite{Wang} established the uniqueness of viscosity solution of $\Delta_\infty$  equation on Carnot groups. Wang \cite{Wang 2} investigated the removable singularities for viscosity subsolutions to degenerate elliptic Pucci operators in the Heisenberg group. 

In this paper, we study viscosity solutions to a class of mixed operators defined by super-positioning Pucci-Heisenberg maximal operator with the fractional sub-Laplacian. We establish the comparison principle, existence and regularity of solutions to the following problem:
\begin{align}\label{eq 2.1}
	\begin{cases}
		\a\mathcal{M}^+_{\lambda,\Lambda}\big(D_{\mathbbm{H}^N,S}^2u\big)-\b(-\Delta_{\mathbbm{H}^N})^su=f &\text{in } {\Omega},\\
		u=g &\text{in }\mathbbm{H}^N\setminus\Omega,
	\end{cases}	
\end{align}
where $\Omega$ is a bounded domain (open, connected set) in Heisenberg group $\mathbbm{H}^N\simeq \mathbb{R}^{2N+1}$, $f\in C(\overline{\Omega}$)$,\,g\in C(\mathbbm{H}^N\setminus\Omega)$ are bounded and  $\a,\b\geq0$ are constants. Here,
\begin{align}\label{mathcal L}
	\mathcal{L}u\coloneqq\a\mathcal{M}^+_{\lambda,\Lambda}\big(D_{\mathbbm{H}^N,S}^2u\big)-\b(-\Delta_{\mathbbm{H}^N})^su	
\end{align}
is a mixed operator on $\mathbbm{H}^N,$ where $\mathcal{M}^+_{\lambda,\Lambda}$ is the Pucci's extremal (maximal) operator which is the fully nonlinear operator and $(-\Delta_{\mathbbm{H}^N})^s$ is the fractional sub-Laplacian.

It is worth noticing that differential equations in the Heisenberg group setting occur naturally in the study of magnetic trajectories on nilmanifolds \cite{ferro}, quantum mechanics \cite{quantum}, non-Markovian coupling of Brownian motions \cite{non-markov} and many others. We recall that in the above contexts, the non-Euclidean geometry occurs naturally. For more insights on pseudodifferential operators in $\mathbbm{H}^N$, we refer to Branson et al. \cite{Branson}, which discuss the Moser-Trudinger inequalities on Cauchy-Riemann (CR) sphere.

We mention that the behaviour of a charged particle in the presence of a \textit{Lorentz force} is governed by the equation of the form$\colon$
\begin{align}
	\nabla_{\gamma'}\gamma'=qF\gamma',
\end{align}
where $\gamma$ is a curve on a Riemannian manifold $(M,g),$ $\nabla$ is known as corresponding Levi-Civita connection and $F$ is some skew symmetric tensor, see \cite{ferro} for the details.

 We elaborate here an interesting application of PDEs in the Heisenberg group. It is well known that the generator of the Brownian motion on $\mathbbm{H}^1$ is the sub-Laplace operator on $\mathbbm{H}^1.$ We highlight the importance of a non-Markovian coupling of Brownian motions in the Heisenberg group, $\mathbbm{H}^1.$ This coupling helps in proving the gradient estimate for harmonic functions on the Heisenberg group, i.e. the solution to the equation
\begin{align}\label{Hei1}
	\Delta_{\mathbbm{H}^1}u=0 \text{ in }\Omega\subset \mathbbm{H}^1,
\end{align}
where $\Delta_{\mathbbm{H}^1}$ is called the sub-Laplace operator. 

\noindent A Markovian coupling of two Markov processes $X$ and $Y$ is a coupling, where for any time instant $t,$ the joint process $\{(X_s,Y_s): s\geq t\}$ conditioned on the filtration $\sigma\{(X_s,Y_s)\}$ is again a coupling of the laws of $X$ and $Y,$ starting from $(X_t,Y_t).$ The coupling of two Markov processes $X$ and $Y$ is said to be \textit{successful} if these couple in a finite time almost surely, i.e.,
\begin{align}
	\tau(X,Y)\coloneqq\inf\{t\geq 0: X_s=Y_s \text{ for all }s\geq t\}\text{ is finite a.s.}
\end{align}
Let us denote the \textit{total variation norm} of any measure $\nu,$ by $\|\nu\|_{\text{TV}}.$ This norm is given by
\begin{align}
	\|\nu\|_{\text{TV}}=\sup\{|\nu(A)|\colon A\text{ is measurable}\}.
\end{align}
Before we get to the application, let us first mention a result due to S. Banerjee and W. S. Kendall \cite{Kend}, which tells that there exists a \textit{maximal} coupling that is again Markovian provided some constraints on the generator of the Markov process and its state space. Also, we recall that in \cite{Kend 2}, authors described an example that for Kolmogorov diffusions defined as a two-dimensional diffusion given by a standard Brownian motion along with its running time integral, given any Markovian coupling, the probability of failing to couple by time $t,$ can not have same order of decay as the total variation distance. They further showed that if the driving Brownian motions start from the same point, then for any Markovian coupling, the coupling rate is at best of order $t^{-1/2}.$ On the other hand, the total variation distance between the corresponding Kolmogorov diffusions decays like $t^{-3/2}.$

\noindent This raises an interesting question in this thread$\colon$ when can we produce non-Markovian couplings that are good enough in the sense that the total-variation distance between $X$ and $Y$ be minimized? The answer to this question comes from \cite{non-markov}. Authors constructed an explicit successful non-Markovian coupling of two copies of the Markov process, which is Brownian motion on the Heisenberg group whose generator is the sub-Laplacian on the Heisenberg group. Moreover, authors use this coupling to produce gradient estimates for harmonic functions on the Heisenberg group, i.e., solutions to \eqref{Hei1}.

It is easy to see that when $\a=0,$ the operator in $\eqref{mathcal L}$ is the fractional sub-Laplacian. Very recently, Palatucci and Piccinini \cite{Palatucci} considered a large class of nonlinear integro-differential operators in the Heisenberg group $\mathbbm{H}^N.$ They proved general Harnack inequalities for solutions of  Dirichlet problem concerning nonlinear integro-differential operators.

 Ferrari and Vecchi \cite{Ferrai} established the H\"older regularity of uniformly continuous and bounded viscosity solutions of degenerate fully nonlinear equations in $\mathbbm{H}^1$. We also refer to \cite{Cutri} for the existence results and Liouville and Harnack type qualitative properties of fundamental solutions. Li and Wang \cite{Wang 2} established a form of comparison principle for sub-elliptic equations in the Heisenberg setting. We refer to \cite{bard1}, where the authors prove a strong comparison principle for degenerate elliptic equations such as those equations which involve the Pucci's extremal operators over Hörmander vector fields. One may also see Manfredini \cite{Manfredi} et al. for the H\"older continuity and boundedness estimates for nonlinear fractional equations in $\mathbbm{H}^N,$ where authors considered equations driven by integro-differential operators whose model is the fractional $p$-Laplacian on Heisenberg group given by
\begin{align}
	\mathcal{L}u(\xi)\coloneqq P.V.\int_{\mathbbm{H}^N}\frac{\big|u(\xi)-u(\eta)\big|^{p-2}(u(\xi)-u(\eta))}{d_0(\eta^{-1}\,o\,\xi)^{Q+sp}}d\eta, \,\,\xi\in\mathbbm{H}^N.
\end{align}
Here, the symbol $P.V.$ in the expression stands for \enquote{in the principal value sense} and $d_0$ denotes a homogeneous norm on $\mathbbm{H}^N$(see Definition \ref{homo}). $Q=2N+2$ is the homogeneous dimension of $\mathbbm{H}^N.$

Motivated by the above works on the mixed operators and recent works on non-Euclidean settings, we consider a class of mixed operators \eqref{mathcal L} on $\mathbbm{H}^N.$  To the best of our knowledge, the mixed operators of type \eqref{mathcal L} have not been considered in the Heisenberg group setting yet.

In order to prove the main theorems of this paper, we first prove the following result:
\begin{thm}[\textbf{{Stability}}]\label{stability}
	For each $n\in\mathbb{N},$ let $u_n\in$ LSC$(\Omega)$ be bounded in $\mathbbm{H}^N$ and satisfy
	\begin{align}
		\a\mathcal{M}^+_{\lambda,\Lambda}\big(D^2_{\mathbbm{H}^N,S}u_n\big)-\b(-\Delta_{\mathbbm{H}^N})^su_n \leq f_n \text{ in }\Omega
	\end{align}
	in the viscosity sense. Let the following be hold:
	\begin{enumerate}
	\item[(i)] $u_n$ converges to $u$ in the $\Gamma$-sense in $\Omega,$
	\item[(ii)] $u_n$ converges to $u$ a.e. in $\mathbbm{H}^N,$
	\item[(iii)] $f_n\longrightarrow f$ locally uniformly in $\Omega.$
	\end{enumerate}	
Then
	\begin{align}
		\a\mathcal{M}^+_{\lambda,\Lambda}\big(D^2_{\mathbbm{H}^N,S}u\big)-\b(-\Delta_{\mathbbm{H}^N})^su \leq f \text{ in }\Omega.
	\end{align}
\end{thm}
Next, we state the comparison principle for viscosity solutions of PDE concerning operators given by \eqref{mathcal L}. The proof is immediate by making use of Lemma \ref{Lemma 1.6} and Lemma \ref{comp} (see, next).
\begin{thm}[\textbf{Comparison Principle}] \label{Comparison}
Let $\Omega$ be a bounded domain in $\mathbbm{H}^N.$ Let $f\in C(\Omega)$ and  $u,\,v$ be bounded USC and LSC functions in $\Omega,$ respectively, which satisfy
\begin{align}
	\a\mathcal{M}^+_{\lambda,\Lambda}\big(D^2_{\mathbbm{H}^N,S}u\big)-\b(-\Delta_{\mathbbm{H}^N})^su \geq f 
\end{align}
and 
\begin{align}
	\a\mathcal{M}^+_{\lambda,\Lambda}\big(D^2_{\mathbbm{H}^N,S}v\big)-\b(-\Delta_{\mathbbm{H}^N})^sv \leq f 
\end{align}
in the viscosity sense in $\Omega.$ Also, if $u\leq v$ in $\mathbbm{H}^N\setminus\Omega$. Then, $u\leq v$ in $\mathbbm{H}^N.$ 
\end{thm}

Further, using comparison principle for sub and super-solutions and then following  Perron's method, we have the following existence result. 
\begin{thm}\label{Theorem 1.3}
	Let $\Omega\subset\mathbbm{H}^N$ be satisfy the exterior Heisenberg ball condition. Let $g\in C(\mathbbm{H}^N\setminus\Omega)$ be bounded in $\mathbbm{H}^N.$ Then, there exists a viscosity solution $u\in C(\overline{\Omega})$ of
	\begin{align}
		\begin{cases}
			\a\mathcal{M}^+_{\lambda,\Lambda}\big(D^2_{\mathbbm{H}^N,S}u\big)-\b(-\Delta_{\mathbbm{H}^N})^su = 0 &\text{ in }\Omega\\
			u=g &\text{ in }\mathbbm{H}^N\setminus\Omega.
		\end{cases}
	\end{align}  
\end{thm}
As mentioned earlier, several authors have established the existence of solutions to Dirichlet problems using the Perron's method. For the Heisenberg group setting, we refer to \cite{Ochoa}, where the authors prove the existence of a viscosity solution for a class of linear second order equations in Heisenberg group. 
It is easy to observe that when $\a=1, \b=0$ in \eqref{mathcal L}, then $\mathcal{L}$ reduces to Pucci-Heisenberg operator. As mentioned before, Ferrari and Vecchi \cite{Ferrai} studied the H\"older behaviour of fully nonlinear equations. In particular, we have the following result$\colon$ \begin{cor}[Theorem 1.3 \cite{Ferrai}]
	Let $u\in C(\mathbbm{H}^1)$ be a bounded and uniformly continuous viscosity solution of
	\begin{align}
		\mathcal{M}_{\lambda,\Lambda}^+\big(D^2_{\mathbbm{H}^N,S}u(\xi)\big)-c(\xi)u(\xi)=f(\xi) \text{ in }\mathbbm{H}^1.
	\end{align}
Let $L_1,$ $L_2,$ $\gamma_1,$ and $\gamma_2$ be positive constants s. t. $\gamma_i\in(0,1],\,i=1,2$ and for any $ \xi,\,\eta\in \mathbbm{H}^1,$
\begin{align}
	\big|c(\xi)-c(\eta)\big|\leq L_1\big|\xi\,o\,\eta^{-1}\big|^{\gamma_1}_{\mathbbm{H}^1},\qquad \big|f(\xi)-f(\eta)\big|\leq L_2\big|\xi\,o\,\eta^{-1}\big|^{\gamma_2}_{\mathbbm{H}^1}.
\end{align} 
Let $c(\xi)$ be positive for all $\xi\in\mathbbm{H}^1$ and
\begin{align}
	\displaystyle{\inf_{{\xi\in B_R^{\mathbbm{H}^1}}(P)}}c(\xi)\coloneqq c_0>0.
\end{align}
Then there exists $\gamma\coloneqq\gamma(c_0,p,L_1,L_2,\Lambda)\in (0,1], \gamma\leq \min\{\gamma_1,\gamma_2\}$ such that
\begin{align}
	\big|u(\xi)-u(\eta)\big|\leq L\big|\xi\,o\,\eta^{-1}\big|_{\mathbbm{H}^1}^\gamma,\,\text{ for }\xi\in\mathbbm{H}^1,
\end{align}
for some $L=L(c_0,P,L_1,L_2,\Lambda)>0.$
\end{cor}

For more such results on regularity in the Euclidean setting, when there is no local term and integro-differential operators are of fractional Laplacian type, we refer to \cite{Silvestre, Lara}. We also mention a very recent work by Biagi et al. \cite{Biagi2}, where authors established the interior Sobolev regularity as well as boundary regularity of Lipschitz type for mixed local and nonlocal operators in the Euclidean framework. More precisely, authors proved the following interior regularity result:
\begin{cor}[Theorem 1.4 \cite{Biagi2}]
	Let $\Omega\subset\mathbb{R}^N$ be a bounded $C^1$ domain. Let $f\in H^k(\Omega)\big(W^{2,k}(\Omega)\big)$ for some integer $k\geq 0.$ Let $u\in H^1(\mathbb{R}^N)$ be a weak solution of 
	\begin{align}
		-\Delta u+(-\Delta)^s u=f \text{ in }\Omega.
	\end{align}
Then $u\in H^{k+2}_{\text{loc}}(\Omega)$.
\end{cor}
Motivated by the above results, we establish the following interior H\"older regularity result.

\begin{thm}[\textbf{Interior regularity}]\label{Holder}
	Let $u$ be a bounded function in $\mathbbm{H}^N$ which satisfy 
	\begin{align}
		\a\mathcal{M}_{\lambda,\Lambda}^+\big(D^2_{\mathbbm{H}^N,S}u\big)-\b(-\Delta_{\mathbbm{H}^N})^su=0 \text{ in }B_{1}^{\mathbbm{H}^N}
	\end{align}
in the viscosity sense. 
Then there exist constants $C$ and $\c=\c(\lambda,\Lambda,N)\in (0,1)$ such that
\begin{align}
	\big|u(\xi)-u(0)\big|\leq C|\xi|_{\mathbbm{H}^N}^\c\|u\|_{\infty,\mathbbm{H}^N},\,\forall\xi\in B_{\frac{1}{2}}^{\mathbbm{H}^N},
\end{align} 
where $B_{r}^{\mathbbm{H}^N}$ denotes the ball of radius $r$ with center at the origin.
\end{thm}

\begin{rem}
	In the Euclidean setting, when there is no local term and the kernel of non local term is of \enquote{fractional Laplacian type}, Caffarelli and Silvestre \cite{Silvestre} studied the regularity for nonlinear integro differential equation with symmetric kernels. Further, Chang-Lara and D\'avila \cite{Lara} established the interior $C^{1, \alpha}$ regularity for nonlinear integro differential equation with non-symmetric kernels. For equations involving more general class of kernels, Silvestre \cite{Silvestre 2}
	 established the H\"older regularity of solutions to equations driven by integro differential operators with non-negative kernel $K$ satisfying $K(x,y)=K(x,-y)$ along with the condition
	\begin{align}
		\text{sup}_x\int_{\mathbb{R}^N}\big(|y|^2\land 1\big)K(x,y)dy<\infty,
	\end{align}
	see Theorem 5.1 \cite{Silvestre 2}.
\end{rem}

The organization of this paper is as follows. In Section 2, we list the basic definitions and introduce the framework in which we work. Section 3 is dedicated to the proofs of our main results.
\section{Preliminaries}
We first recall the briefs about the Heisenberg group $\mathbbm{H}^N.$ The points in $\mathbbm{H}^N$ are denoted by
\begin{align}
	\xi\coloneqq(z,t)=(x_1,\dots,x_N,y_1,\dots,y_N,t)	
\end{align}
and the group $\mathbbm{H}^N$ is defined as the triplet $\big(\mathbb{R}^{N+1},o,\{\Phi_\lambda\}\big)$, where the group law $o$ is defined as follows:
\begin{align}
	\xi\, o\,\xi'&=\big(x+x',y+y',t+t'+2\langle y,x'\rangle- 2\langle x,y'\rangle\big)\\
	&=\bigg(x_1+x_1',\dots,x_N+x_N',y_1+y_1',\dots,y_N+y_N',t+t'+2\sum_{i=1}^{N}\big(y_ix_i'-x_iy_i'\big)\bigg),
\end{align}
where $\langle.,\,.\rangle$ denotes the standard inner product in $\mathbb{R}^N$.
$(\mathbb{R}^{2N+1},o)$ is a Lie group with identity element the origin $\textbf{0}$ and inverse $\xi^{-1}=-\xi$. The dilation group $\{\Phi_\lambda\}_{\lambda>0}$ is given by
\begin{align}
	\Phi(\lambda): \mathbb{R}^{2N+1}\longrightarrow\mathbb{R}^{2N+1}
\end{align}
such that
\begin{align}
	\xi\mapsto\Phi_\lambda(\xi)\coloneqq\big(\lambda x,\lambda y,\lambda^2 t\big).
\end{align} 
$\mathbbm{H}^N$ is also known as Heisenberg-Weyl group in $\mathbb{R}^{2N+1}$. The Jacobian basis of the Heisenberg Lie algebra of $\mathbbm{H}^N$ is given by
\begin{align}
	X_i=\partial_{x_i}+2y_{i}\partial_t, Y_{i}=\partial _{y_i}-2x_i\partial_t, \, 1\leq i\leq N,\, T=\partial_t.
\end{align}
Given a domain $\Omega\subset \mathbbm{H}^N,$ for $u\in C^1(\Omega,\mathbb
{R}),$ the subgradient or the Heisenberg gradient $\nabla_{\mathbbm{H}^N} u$ is defined as follows:
\begin{align}
	\nabla_{\mathbbm{H}^N} u(\xi)\coloneqq\big(X_1u(\xi),\dots,X_{N}u(\xi),Y_1u(\xi),\dots,Y_Nu(\xi)\big).
\end{align}
Also,
\begin{align}
	D^2_{\mathbbm{H}^N,S}u\coloneqq\begin{bmatrix}
		X_1X_1u &\cdots &X_NX_1u &Y_1X_1u &\cdots &Y_NX_1u\\
		\vdots   &\ddots &\vdots &\vdots &\ddots &\vdots \\
		X_1X_Nu &\dots &X_NX_Nu &Y_1X_Nu &\dots &Y_NX_Nu\\
		X_1Y_1u &\dots &X_NY_1u &Y_1Y_1u &\dots &Y_NY_1u\\
		\vdots   &\ddots &\vdots &\vdots &\ddots &\vdots \\
		X_1Y_Nu &\dots &X_NY_Nu &Y_1Y_Nu &\dots &Y_NY_Nu
	\end{bmatrix}_{Sym},
\end{align}
where
\begin{align}
	A_{\text{Sym}}=\frac{1}{2}\big[A+A^T\big],	\text{ for any matrix }A, 
\end{align}
i.e., symmetric part of the matrix $A$.
Now, since
\begin{align}
	[X_i,Y_{i}]&=X_iY_{i}-Y_{i}X_i\\
	&=(\partial_{x_i}+2y_{i}\partial_t)(\partial_{y_{i}}-2x_i\partial_t)-(\partial_{y_{i}}-2x_i\partial_t)(\partial_{x_i}+2y_{i}\partial_t)\\
	&=-4\partial_t,
\end{align}
so it follows that
\begin{align}
	\text{rank}\big(\text{Lie}\{X_1,X_2,\dots,X_{2N},T\}(0,0)\big)=2N+1,
\end{align}
which is the Euclidean dimension of $\mathbbm{H}^N$. We denote by $Q$, the \textit{homogeneous dimension} of $\mathbbm{H}^N$, which is $Q=2N+2.$ The norm on $\mathbbm{H}^N$ is defined by
\begin{align}
	|\xi|_{\mathbbm{H}^N}\coloneqq\bigg[\bigg(\sum_{i=1}^N \big(x_i^2+y_i^2\big)^2\bigg)+t^2\bigg]^{\frac{1}{4}}.
\end{align}
The corresponding distance on $\mathbbm{H}^N$ is defined as follows:
\begin{align}
	d_{\mathbbm{H}^N}(\xi,\hat{\xi})\coloneqq|\hat{\xi}^{-1}o\,\xi|_{\mathbbm{H}^N},
\end{align}
where $\hat{\xi}^{-1}$ is the inverse of $\hat{\xi}$ with respect to to $o$, i.e., $\hat{\xi}^{-1}=-\hat{\xi}$. 

The sub-Laplacian or the Heisenberg Laplacian (also known as Laplacian-Kohn operator), $\Delta_{\mathbbm{H}^N}$ is the self-adjoint operator defined as
\begin{align}
	\Delta_{\mathbbm{H}^N}&\coloneqq\sum_{i=1}^NX_i^2+Y_i^2\\
	&=\sum_{i=1}^N\frac{\partial^2}{\partial x_i^2}+\frac{\partial^2}{\partial y_i^2}+4y_i\frac{\partial^2}{\partial x_i\partial t}-4x_i\frac{\partial^2}{\partial y_i\partial t}+4\big(x_i^2+y_i^2\big)\frac{\partial^2}{\partial t^2}.	
\end{align}
It is useful to observe that
\begin{align}
	\Delta_{\mathbbm{H}^N}=\text{div}\big(\sigma^T\sigma\nabla u\big),
\end{align}
where
\begin{align}
	\sigma=
	\begin{bmatrix}
		I_N &0 &2y\\
		0 &I_N &-2x
	\end{bmatrix}
\end{align}
and $\sigma^T$ is its transpose. Note that
\begin{align}
	A=\sigma^T\sigma=\begin{bmatrix}
		I_N &0 &2y\\
		0 &I_N &-2x\\
		2y &-2x &4\big(|x|^2+|y|^2\big)
	\end{bmatrix}
\end{align}
is a positive semi-definite matrix with det$(A)=0,\,\, \forall \xi\in \mathbbm{H}^N$.

Let us recall the definition of \textit{viscosity sub/super-solution} of (\ref{eq 2.1}) by evaluating the operators in $C^2$ test function $\phi$ $(\psi)$ touching $u$ locally from above (below) and then the final test function $v$ $(w)$ is defined by taking $v=\phi$ $(v=\psi)$ in a small ball and $v=u$ outside. One may see \cite{Barles, Silvestre}, where the notion of viscosity solution for second order elliptic integro-differential equations and fully nonlinear integro-differential equations is given in the Euclidean setting, respectively. 
For the analogous definitions in the Heisenberg group setting, we refer to \cite{Cutri}. We give the definition of viscosity solution for the operator under consideration in the Heisenberg group setting which is consistent with that given in the above mentioned articles.
\begin{defn}
	Let $u:\mathbbm{H}^N\longrightarrow\mathbb{R}$ be a upper semicontinuous (USC) function. Then $u$ is called a \textit{viscosity subsolution} of (\ref{eq 2.1}) if for any $\xi\in\Omega$ and $C^2$ function $\varphi:\overline{U}\longrightarrow\mathbb{R},$ for some neighborhood $U$ of $\xi$ in $\Omega$ such that $\varphi(\xi)=u(\xi)$ and $\varphi(\eta)>u(\eta)$ for $\eta\in U\setminus \{\xi\},$ we have 
	\begin{align}\label{Sub}
		\a\mathcal{M}^+_{\lambda,\Lambda}\big(D_{\mathbbm{H}^N,S}^2v(\xi)\big)-\b(-\Delta_{\mathbbm{H}^N})^sv(\xi) \geq f(\xi),
	\end{align}
	where
	\begin{align}
		v\coloneqq
		\begin{cases}
			\varphi &\text{ in }U,\\
			u &\text{ in } \mathbbm{H}^N\setminus U.
		\end{cases}
	\end{align}
	Moreover, we say $u$ satisfies $\a\mathcal{M}^+_{\lambda,\Lambda}(D_{\mathbbm{H}^N,S}^2u)-\b(-\Delta_{\mathbbm{H}^N})^su \geq f$ in $\Omega$ in the viscosity sense.
\end{defn}

\begin{defn}
	Let $u:\mathbbm{H}^N\longrightarrow\mathbb{R}$ be a lower semicontinuous (LSC) function.  Then $u$ called a \textit{viscosity supersolution} of (\ref{eq 2.1}) if for any $\xi\in\Omega$ and $C^2$ function $\psi:\overline{U}\longrightarrow\mathbb{R},$ for some neighborhood $U$ of $\xi$ in $\Omega$ such that $\psi(\xi)=u(\xi)$ and $\psi(\eta)<u(\eta)$ for $\eta\in U \setminus \{\xi\},$ we have
	\begin{align}
		\a\mathcal{M}^+_{\lambda,\Lambda}\big(D_{\mathbbm{H}^N,S}^2w(\xi)\big)-\b(-\Delta_{\mathbbm{H}^N})^sw(\xi)\leq f(\xi),
	\end{align}
  where
	\begin{align}
		w\coloneqq
		\begin{cases}
			\psi &\text{ in }U,\\
			u &\text{ in } \mathbbm{H}^N\setminus U.
		\end{cases}
	\end{align}
 Moreover, we say $u$ satisfies $\a\mathcal{M}^+_{\lambda,\Lambda}(D_{\mathbbm{H}^N,S}^2u)-\b(-\Delta_{\mathbbm{H}^N})^su \leq f$ in $\Omega$ in the viscosity sense.
\end{defn}
\begin{defn}\label{Vis sol}
	A continuous function $u$ is said to be a \textit{viscosity solution} of (\ref{eq 2.1}) if it is a subsolution as well as a supersolution of (\ref{eq 2.1}).
\end{defn}

Now, we recall the \textit{exterior Heisenberg ball condition}. For the analogous condition in the Euclidean setting, one may see, for instance \cite{Mou1}.
\begin{defn}{\rm \cite[Definition 4.1]{Cutri}}\label{exterior}
	$\Omega$ is said to satisfy the exterior Heisenberg ball condition if there exists $R>0$ such that for any $\xi\in\partial\Omega$ and $0<r\leq R,$ there exists $\eta_\xi^r\in\Omega^c$ satisfying $\overline{B_r^{\mathbbm{H}^N}(\eta_\xi^r)}\cap\overline{\Omega}=\{\xi\}.$ 
\end{defn} 
Next, we recall the definition of sup (inf)-convolution. The construction of convolutions was done by Jensen et al. \cite{JLS} and further developed on Carnot groups by Wang \cite{Wang}.
\begin{defn}
	For an USC function $u,$ the sup-convolution approximation $u^\varepsilon$ is given by
	\begin{align}
		u^\varepsilon(\xi)=\displaystyle{\sup_{\eta\in \mathbbm{H}^N}}\bigg(u(\eta)-\frac{|\xi\,o\,\eta^{-1}|^4_{\mathbbm{H}^N}}{\mathlarger\varepsilon}\bigg).
	\end{align}
\end{defn}
\begin{defn}
	For an LSC function $u$, the inf-convolution approximation $u_\varepsilon$ is given by
	\begin{align}
		u_\varepsilon(\xi)=\displaystyle{\inf_{\eta\in \mathbbm{H}^N}}\bigg(u(\eta)+\frac{|\xi\,o\,\eta^{-1}|^4_{\mathbbm{H}^N}}{\mathlarger\varepsilon}\bigg).
	\end{align}
\end{defn}
It is easy to see that $u^\varepsilon\geq u$ and $u_\varepsilon\leq u.$  

\begin{defn}
	A sequence of LSC functions, $u_k$ is said to $\Gamma$-converge to $u$ in a set $\Omega$ if the following hold:
	\begin{enumerate}
		\item[(i)] For every sequence $\xi_n\longrightarrow \xi$ in $\Omega,$ we have  
	\begin{align}
		\displaystyle{\liminf_{n\longrightarrow \infty}u_n(\xi_n)}\geq u(\xi).
	\end{align}
	\item[(ii)] For every $\xi\in\Omega,$ there exists a sequence $\{\xi_n\}$ converging to $\xi$ in $\Omega$ (known as $\Gamma$-realising sequence) such that
	\begin{align}
		\displaystyle{\limsup_{n\longrightarrow\infty}u_n(\xi_n)}=u(\xi).
	\end{align}
	\end{enumerate}
\end{defn}
This is known as $\Gamma$-limit in literature, see \cite{Silvestre, Progress}. Note that uniform convergence $\implies$ $\Gamma$-convergence. Also, an important property of $\Gamma$-limits is that if $u_n$ converges to $u,$ which has a strict local minimum at $\xi,$ then $u_n$ would have local minimum at $\xi_n$ for a sequence $\xi_n\longrightarrow \xi.$

\begin{defn}\label{homo}{\rm \cite[Definition 2.1]{Palatucci}}
	A homogeneous norm on $\mathbbm{H}^N$ is a continuous function (with respect to Euclidean topology)
	\begin{align}
		d_0:\mathbbm{H}^N\longrightarrow [0,\infty) \text{ s.t. }
	\end{align}
(i) $d_0(\Phi_\lambda(\xi))=\lambda d_0(\xi),\,\forall\lambda>0\, \text{ and } \,\xi\in \mathbbm{H}^N$\\
(ii) $d_0(\xi)=0 \text{ iff }\xi=0.$
\end{defn}
\noindent Moreover, we say that the homogeneous norm is symmetric if
\begin{align}
	d_0(\xi^{-1})=d_0(\xi), \,\forall\xi\in \mathbbm{H}^N.
\end{align}
Throughout the paper, we consider the standard homogeneous norm on $\mathbbm{H}^N.$ For fixed $\xi_0\in \mathbbm{H}^N\text{ and }\, R>0,$
\begin{align}
	B_R^{\mathbbm{H}^N}(\xi_0)\coloneqq\bigg\{\xi\in \mathbbm{H}^N\,;\,|\xi_0^{-1}\,o\,\xi|_{\mathbbm{H}^N}<R\bigg\}
\end{align}
denotes the ball of radius $R$ around $\xi_0$ in $\mathbbm{H}^N.$ One can see that the Jacobian determinant of the dilation $\Phi_\lambda$ is $\lambda^Q$, where $Q=2N+2$, which is the homogeneous dimension of the Heisenberg group. Let us consider $2N\times(2N+1)$ matrix whose rows are the coefficients of the vector fields $X_i$, i.e.,
\begin{align}
	\sigma=\begin{bmatrix}
		I_N &0 &2y\\
		0 &I_N &-2x
	\end{bmatrix}
\end{align}
for $x=(\xi_1,\dots,\xi_N)^T$ and $y=(\xi_{N+1},\dots,\xi_{2N})^T$. Then the Heisenberg gradient of a function $\Phi:\mathbb{R}^{2N+1} \longrightarrow \mathbb{R}$ is given by
\begin{align}
	\nabla_{\mathbbm{H}^N}\Phi=(X_1\Phi,\dots,X_{2N}\Phi)=\sigma(\xi)\nabla\Phi,
\end{align}
where $\nabla$ denotes the usual gradient. Also, the Heisenberg Hessian of $\Phi$ is given by
\begin{align}\label{sigma T}
	D^2_{\mathbbm{H}_N,S}\Phi&=(X_iX_j\Phi)_{Sym}\\
	&=\sigma D^2\Phi\sigma^T,
\end{align}
where \textit{Sym} denotes the symmetrized matrix.
We recall that the fractional sub-Laplacian operator is defined as
\begin{align}
	(-\Delta_{\mathbbm{H}^N})^su(\xi)=-\frac{1}{2}c(N,s)\int_{\mathbbm{H}^N}\frac{u(\xi\,o\,\eta)+u(\xi\,o\,\eta^{-1})-2u(\xi)}{|\eta|^{Q+2s}_{\mathbbm{H}^N}}d\eta,\,\,u\in H^{s}(\mathbbm{H}^N),\,\,\,\xi\in \mathbbm{H}^N, Q=2N+2,
\end{align} 
see 5.1 \cite{Palatucci}. Here, $c(N,s)$ is a positive constant depending on $N$\,\&\,$s.$ For given two parameters $0<\lambda\leq\Lambda,$ Pucci-Heisenberg operators are defined by the composition of Pucci's extremal operators $\mathcal{M}_{\lambda,\Lambda}^{\pm}$ (see 2.2 \cite{Cabre}) with the Heisenberg Hessian
 as follows:
\begin{align}\label{M+}
	\mathcal{M}_{\lambda,\Lambda}^+\big(D^2_{\mathbbm{H}^N,S}u\big)\coloneqq\Lambda\displaystyle{\sum_{e_i\geq 0}e_i}+\lambda\displaystyle{\sum_{e_i< 0}e_i}
\end{align}
\&
\begin{align}
	\mathcal{M}_{\lambda,\Lambda}^-\big(D^2_{\mathbbm{H}^N,S}u\big)\coloneqq\Lambda\displaystyle{\sum_{e_i\leq 0}e_i}+\lambda\displaystyle{\sum_{e_i> 0}e_i},
\end{align}
where $\{e_i\}_{i=1}^{2N}$ are the eigenvalues of the symmetrized horizontal Hessian matrix $D^2_{\mathbbm{H}^N,S}u.$ Let $\mathcal{S}^{2N}$ be denote the set of all real symmetric $2N\times 2N$ matrices. Consider a subset $\mathcal{S}_{\lambda,\Lambda}^{2N}$ of $\mathcal{S}^{2N}$ whose eigenvalues are in $[\lambda,\Lambda].$ These operators are also defined as
\begin{align}\label{max}
\mathcal{M}_{\lambda,\Lambda}^+\big(D^2_{\mathbbm{H}^N,S}u\big)\coloneqq\displaystyle{\max_{M\in \mathcal{S}_{\lambda,\Lambda}^N}}tr\big(MD^2_{\mathbbm{H}^N,S}u\big)
\end{align}
\&
\begin{align}\label{min}
\mathcal{M}_{\lambda,\Lambda}^-\big(D^2_{\mathbbm{H}^N,S}u\big)\coloneqq\displaystyle{\min_{M\in \mathcal{S}_{\lambda,\Lambda}^N}}tr\big(MD^2_{\mathbbm{H}^N,S}u\big).
\end{align}
Furthermore, consider a set $K$ consisting of $2N\times 2N$ matrices $\gamma$ such that $\gamma\gamma^T\in \mathcal{S}^{2N}_{\lambda,\Lambda}.$ Using this, we can re-write \eqref{max} as 
\begin{align}
	\mathcal{M}_{\lambda,\Lambda}^+\big(D^2_{\mathbbm{H}^N,S}u\big)=\displaystyle{\max_{\gamma\in K}}\,\langle \gamma\gamma^T, D^2_{\mathbbm{H}^N,S}u\rangle.
\end{align} 
Now, for any fixed $\gamma\in K,$ we have the following linear operator, say, $L_\gamma$ given by
\begin{align}\label{L beta}
	L_\gamma\big(D^2_{\mathbbm{H}^N,S}u\big)\coloneqq\langle \gamma\gamma^T, D^2_{\mathbbm{H}^N,S}u\rangle.
\end{align}
It is easy to observe that when $\lambda=\Lambda=1,$ the above mentioned operators reduce to the Heisenberg Laplacian $\Delta_{\mathbbm{H}^N,S}.$ 
%

\section{Proofs of main results}
\noindent \textbf{Proof of Theorem \ref{stability}.} Consider a $C^2$ function $\psi$ that touches $u$ from below at $\xi$ in a neighbourhood $U$ in $\Omega.$ By definition of $\Gamma$-convergence, there exists a sequence $\xi_n\longrightarrow \xi$ such that
\begin{align}
\big(u_n-\psi\big)(\xi_n)=\displaystyle{\inf_U}\big(u_n-\psi\big)=\delta_n.	
\end{align}
Therefore, $\delta_n\longrightarrow 0$ as $n\longrightarrow\infty$ and $\psi+\delta_n$ touches $u_n$ at $\xi_n.$ Now, since
	\begin{align}
		\a\mathcal{M}^+_{\lambda,\Lambda}\big(D^2_{\mathbbm{H}^N,S}u_n\big)-\b(-\Delta_{\mathbbm{H}^N})^su_n \leq f_n \text{ in }\Omega,
	\end{align}
so for 
\begin{align}
	w_n\coloneqq\begin{cases}
		\psi+\delta_n &\text{ in }U,\\
		u_n&\text{ in }\mathbbm{H}^N\setminus U,
	\end{cases}
\end{align}
we have by definition of viscosity solution that
\begin{align}
	\a\mathcal{M}^+_{\lambda,\Lambda}\big(D^2_{\mathbbm{H}^N,S}u_n(\xi_n)\big)-\b(-\Delta_{\mathbbm{H}^N})^su_n(\xi_n) \leq f_n(\xi_n) \text{ in }\Omega.
\end{align}
Further, take $\tau\in U$ such that
\begin{align}
	d\big(\tau,\partial U\big)=\displaystyle{\inf_{\widetilde{\tau}\in \partial U}{{\big|\widetilde{\tau}}^{-1}\,o\,\tau\big|_{\mathbbm{H}^N}}}>\rho>0.
\end{align}
This infers
\begin{align}
	\big|&\big(\a\mathcal{M}^+_{\lambda,\Lambda}\big(D^2_{\mathbbm{H}^N,S}w_n(\tau)\big)-\b(-\Delta_{\mathbbm{H}^N})^sw_n(\tau)\big)-\big(\a\mathcal{M}^+_{\lambda,\Lambda}\big(D^2_{\mathbbm{H}^N,S}w(\tau)\big)-\b(-\Delta_{\mathbbm{H}^N})^sw(\tau)\big)\big|\\
	&\qquad= \big|\big(\a\mathcal{M}^+_{\lambda,\Lambda}\big(D^2_{\mathbbm{H}^N,S}w_n(\tau)\big)-\a\mathcal{M}^+_{\lambda,\Lambda}\big(D^2_{\mathbbm{H}^N,S}w(\tau)\big)\big)-\b((-\Delta_{\mathbbm{H}^N})^sw_n(\tau)-(-\Delta_{\mathbbm{H}^N})^sw(\tau))\big|\\
	&\qquad\leq \a \,\displaystyle{\max\bigg(\big|\max_{M\in \mathcal{S}_{\lambda,\Lambda}^N}}tr\big(MD^2_{\mathbbm{H}^N,S}(w_n(\tau)-w(\tau))\big)\big|,\big|\max_{M\in \mathcal{S}_{\lambda,\Lambda}^N}tr\big(MD^2_{\mathbbm{H}^N,S}(w(\tau)-w_n(\tau))\big)\big|\bigg)\\
	&\qquad\qquad\qquad+\b\big|(-\Delta_{\mathbbm{H}^N})^sw_n(\tau)-(-\Delta_{\mathbbm{H}^N})^sw(\tau)\big|\text{ (by \eqref{L beta})}\\
	&\qquad\leq \a \max_{\gamma\in K}\big|L_\gamma\big(D^2(w_n-w)(\tau)\big)\big|+\b\big|(-\Delta_{\mathbbm{H}^N})^s(w_n-w)(\tau)\big|\\
	&\qquad= \a \max_{\gamma\in K}\big|L_\gamma\big(D^2\delta_n(\tau)\big)\big|+\b\big|(-\Delta_{\mathbbm{H}^N})^s(w_n-w)(\tau)\big|.
\end{align}
This further accords
	\begin{align}
\big|&\big(\a\mathcal{M}^+_{\lambda,\Lambda}\big(D^2_{\mathbbm{H}^N,S}w_n(\tau)\big)-\b(-\Delta_{\mathbbm{H}^N})^sw_n(\tau)\big)-\big(\a\mathcal{M}^+_{\lambda,\Lambda}\big(D^2_{\mathbbm{H}^N,S}w(\tau)\big)-\b(-\Delta_{\mathbbm{H}^N})^sw(\tau)\big)\big|\\
&\qquad\leq \frac{\b}{2}c(N,s)\int_{\mathbbm{H}^N}\frac{\big|{(w_n-w)}(\tau\,o\,\eta)+{(w_n-w)}(\tau\,o\,\eta^{-1})-2{(w_n-w)}(\tau)\big|}{|\eta|^{Q+2s}_{\mathbbm{H}^N}}d\eta\\
    &\qquad=\frac{\b}{2}c(N,s)\int_{B_\rho^{\mathbbm{H}^N}}\frac{\big|{(w_n-w)}(\tau\,o\,\eta)+{(w_n-w)}(\tau\,o\,\eta^{-1})-2{(w_n-w)}(\tau)\big|}{|\eta|^{Q+2s}_{\mathbbm{H}^N}}d\eta\\
    &\qquad\qquad\qquad+\frac{\b}{2}c(N,s)\int_{\mathbbm{H}^N\setminus B_\rho^{\mathbbm{H}^N}}\frac{\big|{(w_n-w)}(\tau\,o\,\eta)+{(w_n-w)}(\tau\,o\,\eta^{-1})-2{(w_n-w)}(\tau)\big|}{|\eta|^{Q+2s}_{\mathbbm{H}^N}}d\eta\\
    &\qquad=\frac{\b}{2}c(N,s)\int_{B_\rho^{\mathbbm{H}^N}}\frac{\big|{\delta_n}(\tau\,o\,\eta)+{\delta_n}(\tau\,o\,\eta^{-1})-2{\delta_n}(\tau)\big|}{|\eta|^{Q+2s}_{\mathbbm{H}^N}}d\eta\\
    &\qquad\qquad\qquad+\frac{\b}{2}c(N,s)\int_{\mathbbm{H}^N\setminus B_\rho^{\mathbbm{H}^N}}\frac{\big|{(w_n-w)}(\tau\,o\,\eta)+{(w_n-w)}(\tau\,o\,\eta^{-1})-2{(w_n-w)}(\tau)\big|}{|\eta|^{Q+2s}_{\mathbbm{H}^N}}d\eta\\
    &\qquad=\frac{\b}{2}c(N,s)\int_{\mathbbm{H}^N\setminus B_\rho^{\mathbbm{H}^N}}\frac{\big|{(w_n-w)}(\tau\,o\,\eta)+{(w_n-w)}(\tau\,o\,\eta^{-1})-2{(w_n-w)}(\tau)\big|}{|\eta|^{Q+2s}_{\mathbbm{H}^N}}d\eta.
\end{align}
Now, using the fact that sequence $w_n$ is bounded and that $(w_n-w)(\tau\,o\,\eta)+(w_n-w)(\tau\,o\,\eta^{-1})-(w_n-w)(\tau)\longrightarrow 0$ a.e. along with the fact that $\frac{1}{|\eta|^{Q+2s}}\in L^1(\mathbbm{H}^N\setminus B_\rho),$ we have by the dominated convergence theorem that
\begin{align}\label{eq 3.2}
\left|\big(\a\mathcal{M}^+_{\lambda,\lambda}\big(D^2_{\mathbbm{H}^N,S}w_n(\tau)\big)-\b(-\Delta_{\mathbbm{H}^N})^sw_n(\tau)\big)-\big(\a\mathcal{M}^+_{\lambda,\lambda}\big(D^2_{\mathbbm{H}^N,S}w(\tau)\big)-\b(-\Delta_{\mathbbm{H}^N})^sw(\tau)\big)\right|\longrightarrow 0, 
\end{align}
uniformly in $\tau$ as $n\longrightarrow \infty.$ Equivalently,
\begin{align}
	\a\mathcal{M}^+_{\lambda,\Lambda}\big(D^2_{\mathbbm{H}^N,S}w_n(\tau)\big)-\b(-\Delta_{\mathbbm{H}^N})^sw_n(\tau)\longrightarrow\a\mathcal{M}^+_{\lambda,\Lambda}\big(D^2_{\mathbbm{H}^N,S}w(\tau)\big)-\b(-\Delta_{\mathbbm{H}^N})^sw(\tau),
\end{align}
\text{locally uniformly in} $U.$ Also,
\begin{align}\label{eq 3.3}
	&\left|\big(\a\mathcal{M}^+_{\lambda,\Lambda}\big(D^2_{\mathbbm{H}^N,S}w_n(\xi_n)\big)-\b(-\Delta_{\mathbbm{H}^N})^sw_n(\xi_n)\big)-\big(\a\mathcal{M}^+_{\lambda,\Lambda}\big(D^2_{\mathbbm{H}^N,S}w(\xi)\big)-\b(-\Delta_{\mathbbm{H}^N})^sw(\xi)\big)\right|\\
	&\qquad\leq \left|\big(\a\mathcal{M}^+_{\lambda,\Lambda}\big(D^2_{\mathbbm{H}^N,S}w_n(\xi_n)\big)- \b(-\Delta_{\mathbbm{H}^N})^sw_n(\xi_n)\big)-\big(\a\mathcal{M}^+_{\lambda,\Lambda}\big(D^2_{\mathbbm{H}^N,S}w(\xi_n)\big)-\b(-\Delta_{\mathbbm{H}^N})^sw(\xi_n)\big)\right|\\
	&\qquad\qquad+\left|\big(\a\mathcal{M}^+_{\lambda,\Lambda}\big(D^2_{\mathbbm{H}^N,S}w(\xi_n)\big)- \b(-\Delta_{\mathbbm{H}^N})^sw(\xi_n)\big)-\big(\a\mathcal{M}^+_{\lambda,\Lambda}\big(D^2_{\mathbbm{H}^N,S}w(\xi)\big)- \b(-\Delta_{\mathbbm{H}^N})^sw(\xi)\big)\right|\\
	&\qquad\leq \left|\big(\a\mathcal{M}^+_{\lambda,\Lambda}\big(D^2_{\mathbbm{H}^N,S}w_n(\xi_n)\big)- \b(-\Delta_{\mathbbm{H}^N})^sw_n(\xi_n)\big)-\big(\a\mathcal{M}^+_{\lambda,\Lambda}\big(D^2_{\mathbbm{H}^N,S}w(\xi_n)\big)- \b(-\Delta_{\mathbbm{H}^N})^sw(\xi_n)\big)\right|\\
	&\qquad\qquad+\left|\a\mathcal{M}^+_{\lambda,\Lambda}\big(D^2_{\mathbbm{H}^N,S}w(\xi_n)\big)-\a\mathcal{M}^+_{\lambda,\Lambda}\big(D^2_{\mathbbm{H}^N,S}w(\xi)\big)\big|+\big|\b(-\Delta_{\mathbbm{H}^N})^sw(\xi_n)-\b(-\Delta_{\mathbbm{H}^N})^sw(\xi)\right|.
\end{align}
Clearly, the first term in the R.H.S. of \eqref{eq 3.3} vanishes as $n\longrightarrow\infty$ by \eqref{eq 3.2}. The third term goes to zero as $n\longrightarrow\infty$ by the continuity of $(-\Delta_{\mathbbm{H}^N})w$ in $U.$ 
Also, the second term vanishes by an application of Theorem VI.2.1 \cite{Bhatia}. Thus, we have that
\begin{align}
\a\mathcal{M}^+_{\lambda,\Lambda}\big(D^2_{\mathbbm{H}^N,S}w_n(\xi_n)\big)-\b(-\Delta_{\mathbbm{H}^N}^s)w_n(\xi_n)\longrightarrow\a\mathcal{M}^+_{\lambda,\Lambda}\big(D^2_{\mathbbm{H}^N,S}(w)(\xi)\big)-\b(-\Delta_{\mathbbm{H}^N})^sw(\xi) \text{ as }n\longrightarrow\infty.	
\end{align}
Further, we get
\begin{align}
	\a\mathcal{M}^+_{\lambda,\Lambda}&\big(D^2_{\mathbbm{H}^N,S}w(\xi)\big)-\b(-\Delta_{\mathbbm{H}^N})^sw(\xi)\\
	&\leq \a\mathcal{M}^+_{\lambda,\Lambda}\big(D^2_{\mathbbm{H}^N,S}w_n(\xi_n)\big)-\b(-\Delta_{\mathbbm{H}^N})^sw_n(\xi_n)\\
	&\qquad\qquad+\big|\a\mathcal{M}^+_{\lambda,\Lambda}\big(D^2_{\mathbbm{H}^N,S}w_n(\xi_n)\big)-\b(-\Delta_{\mathbbm{H}^N})^sw_n(\xi_n)-\big(\a\mathcal{M}^+_{\lambda,\Lambda}\big(D^2_{\mathbbm{H}^N,S}w(\xi)\big)-\b(-\Delta_{\mathbbm{H}^N})^sw(\xi)\big)\big|\\
	&\leq f_n(\xi_n)+\big|\a\mathcal{M}^+_{\lambda,\Lambda}\big(D^2_{\mathbbm{H}^N,S}w_n(\xi_n)\big)-\b(-\Delta_{\mathbbm{H}^N})^sw_n(\xi_n)-\big(\a\mathcal{M}^+_{\lambda,\Lambda}\big(D^2_{\mathbbm{H}^N,S}w(\xi)\big)-\b(-\Delta_{\mathbbm{H}^N})^sw(\xi)\big)\big|\\
	&\leq f(\xi)+|f_n(\xi_n)-f(\xi)|\\
	&\qquad\qquad+\big|\a\mathcal{M}^+_{\lambda,\Lambda}\big(D^2_{\mathbbm{H}^N,S}w_n(\xi_n)\big)-\b(-\Delta_{\mathbbm{H}^N})^sw_n(\xi_n)-\big(\a\mathcal{M}^+_{\lambda,\Lambda}\big(D^2_{\mathbbm{H}^N,S}w(\xi)\big)-\b(-\Delta_{\mathbbm{H}^N})^sw(\xi)\big)\big|.
\end{align}
Finally, using \eqref{eq 3.3} together with the fact that $f_n(\xi_n)\longrightarrow f(\xi),$ since $f_n\longrightarrow f$ locally uniformly and $\xi_n\longrightarrow \xi$ proves the claim.\qed

\begin{rem}
	One can see that an analogous result can be obtained for subsolutions using the $\Gamma$-convergence of USC functions. 
\end{rem}

More precisely, we have the following immediate corollary:
\begin{cor}
	For each $n\in\mathbb{N},$ let $u_n\in C(\Omega)$ be bounded in $\mathbbm{H}^N$ and satisfy
	\begin{align}
		\a\mathcal{M}^+_{\lambda,\Lambda}\big(D^2_{\mathbbm{H}^N,S}u_n\big)-\b(-\Delta_{\mathbbm{H}^N})^su_n=f_n \text{ in }\Omega
	\end{align}
	in the viscosity sense. Let the following be hold:
	\begin{enumerate}
		\item[(i)] $u_n$ converges to $u$ in the $\Gamma$-sense in $\Omega,$
	\item[(ii)] $u_n$ converges to $u$ a.e. in $\mathbbm{H}^N,$
	\item[(iii)] $f_n\longrightarrow f$ locally uniformly in $\Omega$.
	\end{enumerate}	
	Then
	\begin{align}
		\a\mathcal{M}^+_{\lambda,\Lambda}\big(D^2_{\mathbbm{H}^N,S}u\big)-\b(-\Delta_{\mathbbm{H}^N})^su=f \text{ in }\Omega.
	\end{align} 
\end{cor}
Next, we prove a lemma, which we use in proving the ensuing results.
\begin{lem}\label{convo}
	Let  $f\in C(\Omega)$ and $u\in \text{USC}(\overline{\Omega})$ be a bounded viscosity subsolution of \eqref{eq 2.1}. Then
	\begin{align}
		\a\mathcal{M}_{\lambda,\Lambda}^+\big(D^2_{\mathbbm{H}^N,S}u^\varepsilon\big)-(-\Delta_{\mathbbm{H}^N})^su^\varepsilon\geq f-d_\mathlarger\varepsilon \text{ in } \Omega_1\Subset\Omega,
	\end{align} 
	where $d_\mathlarger\varepsilon\longrightarrow 0$ in $\Omega_1$ as $\mathlarger{\mathlarger\varepsilon}\longrightarrow 0$ and depends on the modulus of continuity of $f$.
\end{lem}
\noindent \textbf{Proof of Lemma \ref{convo}.}
Let $x_0\in \Omega_1$ and $\varphi$ be a $C^2$ smooth function touching $u^\varepsilon$ from above in some neighbourhood $B^{\mathbbm{H}^N}_r(\xi_0)\subset \Omega_1$ of $\xi_0$. Let us define
\begin{align}
	\Psi\coloneqq\begin{cases}
		\varphi &\text{ in }B^{\mathbbm{H}^N}_r(\xi_0)\\
		u^\varepsilon &\text{ in }\mathbbm{H}^N\setminus B_r^{\mathbbm{H}^N}(\xi_0).
	\end{cases}
\end{align}	
Next, consider a function
\begin{align}\label{Phi}
\Phi(\xi)\coloneqq\Psi\big(\xi_0\,o\,{\xi_0^*}^{-1}\,o\,\xi\big)+\frac{1}{\mathlarger{\mathlarger\varepsilon}}|\xi_0\,o\,{\xi_0^*}^{-1}|^4_{\mathbbm{H}^N},	
\end{align}
where $\xi^*=(x^{1*},x^{2*},\dots,x^{2N*},t^*)\in\Omega$ such that
\begin{align}
	u^\varepsilon(\xi_0)=u(\xi_0^*)-\frac{|\xi_0\,o\,{\xi_0^*}^{-1}|^4_{\mathbbm{H}^N}}{\mathlarger{\mathlarger\varepsilon}}\text{ (see pp. 9 \cite{Li})}.
\end{align}
Now, by definition, we have
\begin{align}\label{def1}
	u^\varepsilon\big(\xi_0\,o\,{\xi_0^*}^{-1}\,o\,\xi\big)&\geq u(\xi)-\frac{|\big(\xi_0\,o\,{\xi_0^*}^{-1}\,o\,\xi\big)\,o\,\xi^{-1}|^4_{\mathbbm{H}^N}}{\mathlarger{\mathlarger\varepsilon}}.
\end{align}
Since
\begin{align}
\big(\xi_0\,o\,{\xi_0^*}^{-1}\,o\,\xi\big)\,o\,\xi^{-1}=\xi_0\,o\,{\xi_0^*}^{-1},
\end{align}
so it implies
\begin{align}
   u^\varepsilon\big(\xi_0\,o\,{\xi_0^*}^{-1}\,o\,\xi\big)&\geq u(\xi)-\frac{|\xi_0\,o\,{\xi_0^*}^{-1}|^4_{\mathbbm{H}^N}}{\mathlarger{\mathlarger\varepsilon}}.
\end{align}
In other words,
\begin{align}
	u(\xi)\leq u^\varepsilon\big(\xi_0\,o\,{\xi_0^*}^{-1}\,o\,\xi\big)+\frac{|\xi_0\,o\,{\xi_0^*}^{-1}|^4_{\mathbbm{H}^N}}{\mathlarger{\mathlarger\varepsilon}}.
\end{align}
Let $\xi$ be close enough to $\xi_0^*,$ then by definition,
\begin{align}
	u(\xi)&\leq \Psi\big(\xi_0\,o\,{\xi_0^*}^{-1}\,o\,\xi\big)+\frac{|\xi_0\,o\,{\xi_0^*}^{-1}|^4_{\mathbbm{H}^N}}{\mathlarger{\mathlarger\varepsilon}}\\
	&=\Phi(\xi)
\end{align}
and $u(\xi_0^*)=\Phi(\xi_0^*).$ It implies that $\Phi(\xi)$ touches $u$ from above at $\xi_0^*$ in $B_r^{\mathbbm{H}^N}(\xi_0^*)$. Next, we consider a function
\begin{align}
	v\coloneqq\begin{cases}
		\Phi &\text{ in }B_r^{\mathbbm{H}^N}(\xi_0^*)\\
		u &\text{ in }\mathbbm{H}^N\setminus B_r^{\mathbbm{H}^N}(\xi_0^*).
	\end{cases}
\end{align}
Then by definition, we have
\begin{align}
	\a\mathcal{M}_{\lambda,\Lambda}^+\big(D^2_{\mathbbm{H}^N,S}v(\xi_0^*)\big)-\b(-\Delta_{\mathbbm{H}^N})^sv(\xi_0^*)\geq f(\xi_0^*).
\end{align}
In other words,
\begin{align}\label{eq 3.6}
	\a\mathcal{M}_{\lambda,\Lambda}^+\big(D^2_{H^N,s}v(\xi_0^*)\big)+\frac{\b}{2}c(N,s)\int_{\mathbbm{H}^N}\frac{v(\xi_0^*\,o\,\eta)+v(\xi_0^*\,o\,\eta^{-1})-2v(\xi_0^*)}{|\eta^{-1}\,o\,\xi|^{Q+2s}_{\mathbbm{H}^N}}d\eta\geq f(\xi_0^*).
\end{align}
It is easy to see that \eqref{Phi} yields
\begin{align}\label{eq 3.7}
	v(\xi^*_0)=\Phi(\xi_0^*)=\Psi(\xi_0)+\frac{1}{\mathlarger{\mathlarger\varepsilon}}|\xi_0\,o\,{\xi_0^*}^{-1}|^4_{\mathbbm{H}^N},
\end{align}
Also,
\begin{align}
	\Phi(\xi_0^*\,o\,\eta)=\Psi(\xi_0\,o\, \eta)+\frac{1}{\mathlarger{\mathlarger\varepsilon}}|\xi_0\,o\,{\xi_0^*}^{-1}|^4_{\mathbbm{H}^N}
\end{align}
and
\begin{align}
	\Phi(\xi_0^*\,o\,\eta^{-1})=\Psi(\xi_0\,o\,\eta^{-1})+\frac{1}{\mathlarger{\mathlarger\varepsilon}}|\xi_0\,o\,{\xi_0^*}^{-1}|^4_{\mathbbm{H}^N}.
\end{align}
It accords
\begin{align}\label{eq 3.8}
	\Phi(\xi_0^*\,o\,\eta)-\Phi(\xi_0^*)=\Psi(\xi_0\,o\, \eta)-\Psi(\xi_0)
\end{align}
\&
\begin{align}
	\Phi(\xi_0^*\,o\,\eta^{-1})-\Phi(\xi_0^*)=\Psi(\xi_0\,o\,\eta^{-1})-\Psi(\xi_0).
\end{align}
Now, we claim that
\begin{align}\label{nabla H}
	\nabla_{\mathbbm{H}^N}\Phi(\xi_0^*)&=\nabla_{\mathbbm{H}^N}\Psi(\xi_0).
\end{align}
Since
\begin{align}
	X_i\Phi(\xi)&=\partial_{x_i}\Phi(\xi)+2y_{i}\partial_t\Phi(\xi)\\
	&=\partial_{x_i}\Phi\big(x^1,\dots,x^{N},y^1,\dots,y^N,t\big)+2y_{i}\partial_t\Phi\big(x^1,\dots,x^{N},y^1,\dots,y^N,t\big)\\
	&=\partial_{x_i}\Psi\big(\xi_0\,o\,{\xi_0^*}^{-1}\,o\,\xi\big)+2y_{i}{\partial_t}\Psi\big(\xi_0\,o\,{\xi_0^*}^{-1}\,o\,\xi\big) \text{ for }1\leq i\leq N,	
\end{align}
so it gives
\begin{align}
	X_i\Phi(\xi_0^{*})=\partial_{x_i}\Psi(\xi_0)+2y_{i}{\partial_t}\Psi(\xi_0)\text{ for }1\leq i\leq N.
\end{align}
Similarly, one may see that
\begin{align}
	Y_i\Phi(\xi_0^*)=\partial_{y_i}\Psi(\xi_0)-2x_i\partial_t\Psi(\xi_0),\,1\leq i\leq N.
\end{align}
Therefore, \eqref{nabla H} holds. Also,
\begin{align}
D^2_{H^N,S}\Phi(\xi_0^*)=D^2_{\mathbbm{H}^N,S}\Psi(\xi_0)	
\end{align}
further implies
\begin{align}\label{eq 3.9}
	\mathcal{M}_{\lambda,\Lambda}^+\big(D^2_{\mathbbm{H}^N,S}\Phi(\xi_0^*)\big)=\mathcal{M}_{\lambda,\Lambda}^+\big(D^2_{\mathbbm{H}^N,S}\Psi(\xi_0)\big)=\mathcal{M}_{\lambda,\Lambda}^+\big(D^2_{\mathbbm{H}^N,S}\varphi(\xi_0)\big).
\end{align}
Using \eqref{eq 3.8} \& \eqref{eq 3.9} in \eqref{eq 3.6} yields
\begin{align}\label{eq 3.10}
	\quad\quad f(\xi_0^*)&\leq \a\mathcal{M}_{\lambda,\Lambda}^+\big(D^2_{\mathbbm{H}^N,S}v(\xi_0^*)\big)+\frac{\b}{2}c(N,s)\int_{\mathbbm{H}^N}\frac{v(\xi_0^*\,o\,\eta)+v(\xi_0^*\,o\,\eta^{-1})-2v(\xi_0^*)}{|\eta|^{Q+2s}_{\mathbbm{H}^N}}d\eta\\
	&=\a\mathcal{M}_{\lambda,\Lambda}^+\big(D^2_{\mathbbm{H}^N,S}\varphi(\xi_0)\big)+\frac{\b}{2}c(N,s)\int_{\{\eta\in \mathbbm{H}^N:\,\xi_0^*\,o\,\eta,\,\xi_0^*\,o\,\eta^{-1}\in B_r^{\mathbbm{H}^N}(\xi_0^*)\}}\frac{\Phi(\xi_0^*\,o\,\eta)+\Phi(\xi_0^*\,o\,\eta^{-1})-2\Phi(\xi_0^*)}{|\eta|^{Q+2s}_{\mathbbm{H}^N}}d\eta\\
	&\quad\qquad\qquad\qquad\qquad\qquad +\frac{\b}{2}c(N,s)\int_{\mathbbm{H}^N\setminus\{\eta\in \mathbbm{H}^N:\,\xi_0^*\,o\,\eta,\,\xi_0^*\,o\,\eta^{-1}\in B_r^{\mathbbm{H}^N}(\xi_0^*)\}}\frac{u(\xi_0^*\,o\,\eta)+u(\xi_0^*\,o\,\eta^{-1})-2u(\xi_0^*)}{|\eta|^{Q+2s}_{\mathbbm{H}^N}}d\eta\\
	&=\a\mathcal{M}_{\lambda,\Lambda}^+\big(D^2_{\mathbbm{H}^N,S}\varphi(\xi_0)\big)+\frac{\b}{2}c(N,s)\int_{\{\eta\in \mathbbm{H}^N:\,\xi_0^*\,o\,\eta,\,\xi_0^*\,o\,\eta^{-1}\in B_r^{\mathbbm{H}^N}(\xi_0^*)\}}\frac{\Phi(\xi_0^*\,o\,\eta)+\Phi(\xi_0^*\,o\,\eta^{-1})-2\Phi(\xi_0^*)}{|\eta|^{Q+2s}_{\mathbbm{H}^N}}d\eta\\
	&\quad\qquad\qquad\qquad\qquad\qquad +\frac{\b}{2}c(N,s)\int_{\mathbbm{H}^N\setminus\{\eta\in \mathbbm{H}^N:\,\xi_0^*\,o\,\eta,\,\xi_0^*\,o\,\eta^{-1}\in B_r^{\mathbbm{H}^N}(\xi_0^*)\}}\frac{u(\xi_0^*\,o\,\eta)+u(\xi_0^*\,o\,\eta^{-1})-2\Phi(\xi_0^*)}{|\eta|^{Q+2s}_{\mathbbm{H}^N}}d\eta.
\end{align}
By the definition of $u^\varepsilon,$ we have
\begin{align}
	u(\eta)\leq u^\varepsilon(\xi)+\frac{\big|\xi\,o\,{\eta}^{-1}\big|^4_{\mathbbm{H}^N}}{\mathlarger{\mathlarger\varepsilon}} \text{ for all }\xi,\,\eta\in \mathbbm{H}^N.
\end{align}
It further implies
\begin{align}\label{use compo}
	u(\xi_0^*\,o\,\eta)\leq u^\varepsilon(\xi_0\,o\,\eta)+\frac{\big|{(\xi_0\,o\,\eta)\,o\,(\xi_0^*\,o\,\eta)}^{-1}\big|^4_{\mathbbm{H}^N}}{\mathlarger{\mathlarger\varepsilon}},
\end{align}
and
\begin{align}\label{use compo 2}
	u(\xi_0^*\,o\,\eta^{-1})\leq u^\varepsilon(\xi_0\,o\,\eta^{-1})+\frac{\big|(\xi_0\,o\,\eta^{-1})\,o\,{(\xi_0^*\,o\,\eta^{-1})}^{-1}\big|^4_{\mathbbm{H}^N}}{\mathlarger{\mathlarger\varepsilon}}.
\end{align}
Let 
\begin{align}
\xi_0^*&=(x_0^*,y_0^*,t_0^*)=\big(x_0^{1*},\dots,x_0^{N*},y_0^{1*},\dots,y_0^{N*},t_0^*\big),\\
\xi_0&=(x_0,y_0,z_0,t)=\big(x^1_0,\dots,x^N_0,y^1_0,\dots,y^N_0,t_0\big),\\
\eta&=(x,y,z,t)=\big(x_1,\dots,x_N,y_1,\dots,y_N,t\big).
\end{align}
It infers that
\begin{align}
	\xi_0\,o\,\eta&=\big(x_0^1+x^1,\dots,x_0^N+x^N,y_0^1+y^1,\dots,y_0^N+y^N,t_0+t+2\langle y_0,x\rangle-2\langle x_0,y\rangle\big),\\
	\xi_0^*\,o\,\eta&=\big(x_0^{1*}+x^{1},\dots,x_0^{N*}+x^{N},y_0^{1*}+y^{1},\dots,y^{N*}_0+y^{N},t_0^*+t+2\langle y_0^*,x\rangle-2\langle x_0^*,y\rangle\big).
\end{align}
This immediately gives
\begin{align}
	{(\xi_0^*\,o\,\eta)}^{-1}=\big(-x_0^{1*}-x^{1},\dots,-x_0^{N*}-x^{N},-y_0^{1*}-y^{1},\dots,-y_0^{N*}-y^{N},-t_0^*-t-2\langle y_0^*,x\rangle+2\langle x_0^*,y\rangle\big),
\end{align}
which further grants 
\begin{align}\label{compo}
	\qquad\qquad(\xi_0\,o\,\eta)\,o\,{(\xi_0^*\,o\,\eta)}^{-1}&=\big(x_0^1-x_0^{1*},\dots,x_0^N-x_0^{N*},y_0^1-y_0^{1*},\dots,y_0^N-y_0^{N*},t_0-t_0^*\\
	&\qquad+2\langle y,x_0^*\rangle-2\langle x,y_0^*\rangle-2\langle y,x_0\rangle+2\langle x,y_0\rangle+2\langle y_0+y,-x_0^*-x\rangle-2\langle x_0+x,-y_0^*-y\rangle\big)\\
	&=\big(x_0^1-x_0^{1*},\dots,x_0^N-x_0^{N*},y_0^1-y_0^{1*},\dots,y_0^N-y_0^{N*},t_0-t_0^*\\
	&\qquad+2\langle y,x_0^*\rangle-2\langle x,y_0^*\rangle-2\langle y,x_0\rangle+2\langle x,y_0\rangle+2\langle y_0,-x_0^*-x\rangle\\
	&\qquad+2\langle y,-x_0^*-x\rangle-2\langle x_0,-y_0^*-y\rangle-2\langle x,-y_0^*-y\rangle\big)\\
	&=\big(x_0^1-x_0^{1*},\dots,x_0^N-x_0^{N*}+,y_0^1-y_0^{1*},\dots,y_0^N-y_0^{N*},t_0-t_0^*\\
	&\qquad+2\langle y,x_0^*\rangle-2\langle x,y_0^*\rangle-2\langle y,x_0\rangle+2\langle x,y_0\rangle+2\langle y_0,-x_0^*\rangle
	+2\langle y_0,-x\rangle+2\langle y,-x_0^*\rangle
	+2\langle y,-x\rangle\\
	&\qquad-2\langle x_0,-y_0^*\rangle-2\langle x_0,-y\rangle-2\langle x,-y_0^*\rangle-2\langle x,-y\rangle\big)\\
	&=\big(x_0^1-x_0^{1*},\dots,x_0^N-x_0^{N*}+,y_0^1-y_0^{1*},\dots,y_0^N-y_0^{N*},t_0-t_0^*\\
	&\qquad+2\langle y,x_0^*\rangle-2\langle x,y_0^*\rangle-2\langle y,x_0\rangle+2\langle x,y_0\rangle+2\langle x_0,y_0^*\rangle+2\langle x_0,y\rangle+2\langle x,y_0^*\rangle+2\langle x,y\rangle\\
	&\qquad-2\langle y_0,x_0^*\rangle
	-2\langle y_0,x\rangle-2\langle y,x_0^*\rangle
	-2\langle y,x\rangle\big)\\
	&=\big(x_0-x_0^*,y_0-y_0^*,t_0-t_0^*+2\langle x_0,y_0^*\rangle-\langle y_0,x_0^* \rangle\big)\\
	&=\big(x_0-x_0^*,y_0-y_0^*,t_0-t_0^*+2\langle y_0,-x_0^*\rangle-2\langle x_0,-y_0^* \rangle\big)\\
	&=\xi_0\,o\,{\xi_0^*}^{-1}.
\end{align}
Similarly, one may see that
\begin{align}\label{compo 2}
	(\xi_0\,o\,\eta^{-1})\,o\,{(\xi_0^*\,o\,\eta^{-1})}^{-1}=\xi_0\,o\,{\xi_0^*}^{-1}.
\end{align}
Using \eqref{compo} \& \eqref{compo 2} in \eqref{use compo} \& \eqref{use compo 2}, respectively infers
\begin{align}\label{eq 3.11}
	\qquad\quad u(\xi_0^*\,o\,\eta)+u(\xi_0^*\,o\,\eta^{-1})-2\Phi(\xi_0^*)&\leq u^\varepsilon(\xi_0\,o\,\eta)+u^\varepsilon(\xi_0\,o\,\eta^{-1})+\frac{2\big|\xi_0\,o\,{\xi_0^*}^{-1}\big|^4_{\mathbbm{H}^N}}{\mathlarger{\mathlarger\varepsilon}}-2\Phi(\xi_0^*)\\
	&=u^\varepsilon(\xi_0\,o\,\eta)+u^\varepsilon(\xi_0\,o\,\eta^{-1})+\frac{2\big|\xi_0\,o\,{\xi_0^*}^{-1}\big|^4_{\mathbbm{H}^N}}{\mathlarger{\mathlarger\varepsilon}}-2\bigg(\Psi(\xi_0)+\frac{1}{\mathlarger{\mathlarger\varepsilon}}\big|\xi_0\,o\,{\xi_0^*}^{-1}\big|^4_{\mathbbm{H}^N}\bigg)\text{ (by \eqref{eq 3.7})}\\
	&=u^\varepsilon(\xi_0\,o\,\eta)+u^\varepsilon(\xi_0\,o\,\eta^{-1})-2\Psi(\xi_0).
\end{align}
We also have
\begin{align}\label{eq 3.12}
	\big|\xi_0\,o\,{\xi_0^*}^{-1}\big|^4_{\mathbbm{H}^N}\leq \mathlarger\varepsilon\,\underset{\Omega}{\text{osc}}~u,
\end{align}
see \cite{Li, Wang} for the details. It immediately grants that $|\xi_0\,o\,{\xi_0^*}^{-1}|_{\mathbbm{H}^N}$ can be made as small as possible by the suitable choice of $\mathlarger{\mathlarger\varepsilon}$. Using \eqref{eq 3.11} in \eqref{eq 3.10}, we get
\begin{align}\label{eq 3.13}
\quad\qquad f(\xi_0^*)&\leq \a\mathcal{M}_{\lambda,\Lambda}^+\big(D^2_{\mathbbm{H}^N,S}\varphi(\xi_0)\big)+\frac{\beta}{2}c(N,s)\int_{\{\eta\in \mathbbm{H}^N:\,\xi_0^*\,o\,\eta,\,\xi_0^*\,o\,\eta^{-1}\in B_r^{\mathbbm{H}^N}(\xi_0^*)\}}\frac{\Phi(\xi_0^*\,o\,\eta)+\Phi(\xi_0^*\,o\,\eta^{-1})-2\Phi(\xi_0^*)}{|\eta|^{Q+2s}_{\mathbbm{H}^N}}d\eta\\
&\quad\qquad\qquad\qquad\qquad\qquad+\frac{\beta}{2}c(N,s)\int_{\mathbbm{H}^N\setminus\{\eta\in \mathbbm{H}^N:\,\xi_0^*\,o\,\eta,\,\xi_0^*\,o\,\eta^{-1}\in B_r^{\mathbbm{H}^N}(\xi_0^*)\}}\frac{u(\xi_0^*\,o\,\eta)+u(\xi_0^*\,o\,\eta^{-1})-2\Phi(\xi_0^*)}{|\eta|^{Q+2s}_{\mathbbm{H}^N}}d\eta\\
&\leq \a\mathcal{M}_{\lambda,\Lambda}^+\big(D^2_{\mathbbm{H}^N,S}\varphi(\xi_0)\big)+c(N,s)\int_{\{\eta\in \mathbbm{H}^N:\,\xi_0^*\,o\,\eta,\,\xi_0^*\,o\,\eta^{-1}\in B_r^{\mathbbm{H}^N}(\xi_0^*)\}}\frac{\varphi(\xi_0\,o\,\eta)+\varphi(\xi_0\,o\,\eta^{-1})-2\varphi(\xi_0)}{|\eta|^{Q+2s}_{\mathbbm{H}^N}}d\eta\\
&\quad\qquad\qquad\qquad\qquad\qquad+\frac{\beta}{2}c(N,s)\int_{\mathbbm{H}^N\setminus\{\eta\in \mathbbm{H}^N:\,\xi_0^*\,o\,\eta,\,\xi_0^*\,o\,\eta^{-1}\in B_r^{\mathbbm{H}^N}(\xi_0^*)\}}\frac{u^\varepsilon(\xi_0\,o\,\eta)+u^\varepsilon(\xi_0\,o\,\eta^{-1})-2\Psi(\xi_0)}{|\eta|^{Q+2s}_{\mathbbm{H}^N}}d\eta\\
&=\a\mathcal{M}_{\lambda,\Lambda}^+\big(D^2_{\mathbbm{H}^N,S}\varphi(\xi_0)\big)-\b(-\Delta_{\mathbbm{H}^N})^s\Psi(\xi_0).
\end{align}
Now, adding and subtracting $f(\xi_0)$ in the L.H.S. of \eqref{eq 3.13} gives
\begin{align}
	\a\mathcal{M}_{\lambda,\Lambda}^+\big(D^2_{\mathbbm{H}^N,S}\Psi(\xi_0)\big)-\b(-\Delta_{\mathbbm{H}^N})^s\Psi(\xi_0)&\geq f(\xi_0)-\big(f(\xi_0)-f(\xi_0^*)\big)\\
	&\geq f(\xi_0)-|f(\xi_0)-f(\xi_0^*)|.
\end{align}
Next, using \eqref{eq 3.12} together with the continuity of $f$, we define
\begin{align}
	d_\varepsilon=\displaystyle{\sup_{B^{\mathbbm{H}^N}_{\delta\sqrt{\varepsilon}}(\xi_0)}}|f(\xi_0)-f(\xi_0^*)|,
\end{align}
for some $\delta=\delta(\underset{\Omega}{\text{osc}}\,u)$. Clearly, $d_\varepsilon\longrightarrow 0$ as $\mathlarger\varepsilon\longrightarrow 0.$ Hence the claim.\qed

\begin{rem}
	Using the similar arguments, one may see that an analogous result also holds for supersolutions.
\end{rem}
\begin{lem}\label{Lemma 1.6}
	Let $f$ $\&$ $g$ be two continuous functions. Let $u$ $\&$ $v$ be bounded USC and $LSC$ functions in $\mathbbm{H}^N$, respectively. Let
	\begin{align}
		\a\mathcal{M}^+_{\lambda,\Lambda}\big(D^2_{\mathbbm{H}^N,S}u\big)-\b(-\Delta_{\mathbbm{H}^N})^su\geq f
	\end{align} 	
	$\&$
	\begin{align}
		\a\mathcal{M}^+_{\lambda,\Lambda}\big(D^2_{\mathbbm{H}^N,S}v\big)-\b(-\Delta_{\mathbbm{H}^N})^sv\leq g
	\end{align}
	be hold in the viscosity sense in $\Omega.$ Then
	\begin{align}
		\a\mathcal{M}^+_{\lambda,\Lambda}\big(D^2_{\mathbbm{H}^N,S}(u-v)\big)-\b(-\Delta_{\mathbbm{H}^N})^s(u-v)\geq f-g \text{ in }\Omega \text{ in the viscosity sense.}
	\end{align}
\end{lem}

\noindent \textbf{Proof of Lemma \ref{Lemma 1.6}.}
By Lemma \ref{convo}, we have
\begin{align}\label{1}
	\a\mathcal{M}_{\lambda,\Lambda}^+\big(D^2_{\mathbbm{H}^N,S}u^\varepsilon\big)-\b(-\Delta_{\mathbbm{H}^N})^su^\varepsilon\geq f-d_\mathlarger\varepsilon	
\end{align}
$\&$ 
\begin{align}\label{2}
	\a\mathcal{M}_{\lambda,\Lambda}^+\big(D^2_{\mathbbm{H}^N,S}v_\varepsilon\big)-\b(-\Delta_{\mathbbm{H}^N})^sv_\varepsilon\leq g+d_\mathlarger\varepsilon.
\end{align}
Our aim is to show that
\begin{align}
	\a\mathcal{M}_{\lambda,\Lambda}^+\big(D^2_{\mathbbm{H}^N,S}(u^\varepsilon-v_\varepsilon)\big)-\b(-\Delta_{\mathbbm{H}^N})^s(u^\varepsilon-v_\varepsilon)\geq f-g-2d_\mathlarger\varepsilon,	
\end{align}
so that further using the stability result (Theorem \ref{stability}) we get the claim. 

\noindent For this, let $P\in C_b^2(\mathbbm{H}^N)$ be a paraboloid in $\overline{B^{\mathbbm{H}^N}_r}\subset \Omega_1$ such that
\begin{align}
	P(\xi_0)=u^\varepsilon(\xi_0)-v_\varepsilon(\xi_0)
\end{align}
and
\begin{align}
	P(\tau)\geq u^\varepsilon(\tau)-v_\varepsilon(\tau) \text{ for all }\tau \in B^{\mathbbm{H}^N}_r(\xi_0).
\end{align}
We assume that $\overline{B^{\mathbbm{H}^N}_{2r}}(\xi_0)\subset\Omega.$ Consider
\begin{align}
\Phi(x)\coloneqq\begin{cases}
	P &\text{ in }B_r^{\mathbbm{H}^N}(\xi_0)\\
	u^\varepsilon-v_\varepsilon &\text{ in }\mathbbm{H}^N\setminus B_r^{\mathbbm{H}^N}(\xi_0).
\end{cases}	
\end{align}
Take $\delta>0$ and let us define
\begin{align}
	w(\xi)=v_\varepsilon(\xi)-u^\varepsilon(\xi)+\Phi(\xi)+\delta\left(\big|{\xi_0}^{-1}o\,\xi\,\big|_{\mathbbm{H}^N}\land r\right)^4-\delta r_1^4,
\end{align}
for $0<r_1<\delta\land\frac{r}{2}.$ We observe that $w\geq 0$ on $\partial B_{r_1}^{\mathbbm{H}^N}(\xi_0)$ and $w(\xi_0)=-\delta r^2<0$. By Theorem 5.1 \cite{Cabre}, we have that for any $\xi\in \overline{B_{r_1}^{\mathbbm{H}^N}}(\xi_0),$ there exists a convex paraboloid $P^\xi$ of opening $K$ (some constant independent of $\xi$) which touches $w$ from above at $\xi \in B_{r_1}^{\mathbbm{H}^N}(\xi_0)$. Further, by Lemma 3.5 \cite{Cabre} and $w(\xi_0)<0,$ we have that 
\begin{align}\label{eq det}
	0<\int_{B_r^{\mathbbm{H}^N}(\xi_0)\cap\{w=\Gamma_w\}}\det D^2\Gamma_w,
\end{align} 
where $\Gamma_w$ is the convex envelope of $w$ in $A$ given by
\begin{align}
	\Gamma_w(\xi)=\displaystyle{\sup_v}\bigg\{v(\xi)\,:\,v\leq w \text{ in }A, v \text{ convex in }A\bigg\},
\end{align}
for $\xi\in A$. Here, $A\subset B_r^{\mathbbm{H}^N}(\xi_0)$ with $|B_r^{\mathbbm{H}^N}(\xi_0)\setminus A|=0.$ Moreover, $u^\varepsilon$ and $v_\varepsilon$ are punctually second order differentiable in $A$ (see 2.2 (ii) \cite{Li}). It gives that $\mathcal{M}_{\lambda,\Lambda}^+(D^2_{\mathbbm{H}^N,S}u^\varepsilon)-(-\Delta_{\mathbbm{H}^N})^su^\varepsilon$ and $\mathcal{M}_{\lambda,\Lambda}^+(D^2_{\mathbbm{H}^N,S}v_\varepsilon)-(-\Delta_{\mathbbm{H}^N})^sv_\varepsilon$ are defined in the classical sense for $\xi\in A.$

Now, using the convexity of $\Gamma_w$ and that $\Gamma_w\leq w,$ we have that the Hessian matrix $D^2w(\xi)$ is semi-positive definite for $\xi\in A\cap\{w=\Gamma_w\}.$ It yields using \eqref{sigma T} that
\begin{align}\label{eq 3.18}
\mathcal{M}_{\lambda,\Lambda}^-\big(D^2_{\mathbbm{H}^N,S}w(\xi)\big)\geq 0,	
\end{align} 
and
\begin{align}\label{eq 3.20}
	\quad \qquad\frac{1}{2}c(N,s)&\int_{\{\eta\in \mathbbm{H}^N:\,\xi\,o\,\eta\in B_r^{\mathbbm{H}^N}(\xi_0)\}}\frac{w(\xi\,o\,\eta)+w(\xi\,o\,\eta^{-1})-2w(\xi)}{|\eta|^{Q+2s}_{\mathbbm{H}^N}}d\eta\\
	&\qquad\qquad\qquad=\frac{1}{2}c(N,s)\int_{\{\eta\in \mathbbm{H}^N:\,\xi\,o\,\eta\in B_r^{\mathbbm{H}^N}(\xi_0)\}}\frac{w(\xi\,o\,\eta)-w(\xi)-\mathbbm{1}_{\{|\eta|_{\mathbbm{H}^N}\leq 1\}}\eta\cdot(\nabla_{\mathbbm{H}^N}w(\xi),\,\partial_tw(\xi))}{|\eta|^{Q+2s}_{\mathbbm{H}^N}}d\eta\\
	&\qquad\qquad\qquad\geq 0,
\end{align}
for $\xi\in A$ using Proposition 2.4 \cite{Palatucci} and
\begin{align}
	w(\xi\,o\,\eta)-w(\xi)-\mathbbm{1}_{\{|\eta|_{\mathbbm{H}^N}\leq 1\}}\eta\cdot (\nabla_{\mathbbm{H}^N}w(\xi),\,\partial_tw(\xi))&\geq \Gamma_w(\xi\,o\,\eta)-\Gamma_w(\xi)-\mathbbm{1}_{\{|\eta|_{\mathbbm{H}^N}\leq 1\}}\eta\cdot(\nabla_{\mathbbm{H}^N}\Gamma_w(\xi),\,\partial_t\Gamma_w(\xi))\\
	&\geq 0
\end{align}
along with $\nabla_ w(\xi)=\nabla\Gamma_w(\xi)$ for $\xi
\in \{w=\Gamma_w\}.$
It is clear from \eqref{eq det} and $\big|B_r^{{\mathbbm{H}^N}}(\xi_0)\setminus A\big|=0$ that
\begin{align}
	\big|\{w=\Gamma_w\}\cap A\big|>0,
\end{align}
i.e., there is a point $\xi_\delta\in\{w=\Gamma_w\}\cap A$ such that \eqref{1} and \eqref{2} hold classically. Let $I_N$ be denote the identity matrix of order $N$. We have
\begin{align}\label{11}
	\\
	f(\xi_\delta)-d_\varepsilon&\leq 	\a\mathcal{M}_{\lambda,\Lambda}^+\big(D^2_{\mathbbm{H}^N,S}u^\varepsilon(\xi_\delta)\big)-\b(-\Delta_{\mathbbm{H}^N})^su^\varepsilon(\xi_\delta)\\
	&= \a\mathcal{M}_{\lambda,\Lambda}^+\big(D^2_{\mathbbm{H}^N,S}v_\varepsilon(\xi_\delta)-D^2_{\mathbbm{H}^N,S}w(\xi_\delta)+D^2_{\mathbbm{H}^N,S}\Phi(\xi_\delta)+D^2_{\mathbbm{H}^N,S}\big(\delta\big|{\xi_0}^{-1}o\,\xi\,\big|^4_{\mathbbm{H}^N}(\xi_\delta)\big)\big)-\b(-\Delta_{\mathbbm{H}^N})^su^\varepsilon(\xi_\delta)\\
	&=\a\mathcal{M}_{\lambda,\Lambda}^+\big(D^2_{\mathbbm{H}^N,S}v_\varepsilon(\xi_\delta)+D^2_{\mathbbm{H}^N,S}\Phi(\xi_\delta)+\delta D^2_{\mathbbm{H}^N,S}\big(\big|{\xi_0}^{-1}o\,\xi\,\big|^4_{\mathbbm{H}^N}(\xi_\delta)\big)\big)-D^2_{\mathbbm{H}^N,S}w(\xi_\delta)\big)-\b(-\Delta_{\mathbbm{H}^N})^su^\varepsilon(\xi_\delta)\\
	&\leq \a\mathcal{M}_{\lambda,\Lambda}^+\big(D^2_{\mathbbm{H}^N,S}v_\varepsilon(\xi_\delta)+D^2_{\mathbbm{H}^N,S}\Phi(\xi_\delta)+\delta D^2_{\mathbbm{H}^N,S}\big(\big|{\xi_0}^{-1}o\,\xi\,\big|^4_{\mathbbm{H}^N}(\xi_\delta)\big)\big)\\
	&\qquad+\a\mathcal{M}_{\lambda,\Lambda}^+\big(-D^2_{\mathbbm{H}^N,S}w(\xi_\delta)\big)-\b(-\Delta_{\mathbbm{H}^N})^su^\varepsilon(\xi_\delta)\\
	&=\a\mathcal{M}_{\lambda,\Lambda}^+\big(D^2_{\mathbbm{H}^N,S}v_\varepsilon(\xi_\delta)+D^2_{\mathbbm{H}^N,S}\Phi(\xi_\delta)+\delta D^2_{\mathbbm{H}^N,S}\big(\big|{\xi_0}^{-1}o\,\xi\,\big|^4_{\mathbbm{H}^N}(\xi_\delta)\big)\big)\\&\qquad-\a\mathcal{M}_{\lambda,\Lambda}^-\big(D^2_{\mathbbm{H}^N,S}w(\xi_\delta)\big)-\b(-\Delta_{\mathbbm{H}^N})^su^\varepsilon(\xi_\delta)\\
	&\leq \a\mathcal{M}_{\lambda,\Lambda}^+\big(D^2_{\mathbbm{H}^N,S}v_\varepsilon(\xi_\delta)+D^2_{\mathbbm{H}^N,S}\Phi(\xi_\delta)+\delta D^2_{\mathbbm{H}^N,S}\big(\big|{\xi_0}^{-1}o\,\xi\,\big|^4_{\mathbbm{H}^N}(\xi_\delta)\big)\big)-\b(-\Delta_{\mathbbm{H}^N})^su^\varepsilon(\xi_\delta)\text{ (by \eqref{eq 3.18})},
\end{align}
where in the second last step, we used the relation
\begin{align}
	\mathcal{M}_{\lambda,\Lambda}^+(-M)=-\mathcal{M}_{\lambda,\Lambda}^-(M).
\end{align}
Now, since $\big|{\xi_0}^{-1}o\,\xi\,\big|^4_{\mathbbm{H}^N}$ is a radial function so using Lemma 3.2 \cite{Cutri} together with sub-additivity property of Pucci's maximal operator and \eqref{11} yields
\begin{align}
	f(\xi_\delta)-d_\varepsilon&\leq \a\mathcal{M}_{\lambda,\Lambda}^+\big(D^2_{\mathbbm{H}^N,S}v_\varepsilon(\xi_\delta)\big)+\a\mathcal{M}_{\lambda,\Lambda}^+\big(D^2_{\mathbbm{H}^N,S}\Phi(\xi_\delta)\big)+\a\delta \mathcal{M}_{\lambda,\Lambda}^+\big(D^2_{\mathbbm{H}^N,S}\big(\big|{\xi_0}^{-1}o\,\xi_\delta\,\big|^4_{\mathbbm{H}^N}\big)\big)-\b(-\Delta_{\mathbbm{H}^N})^su^\varepsilon(\xi_\delta)\\
	&=\a\mathcal{M}_{\lambda,\Lambda}^+\big(D^2_{\mathbbm{H}^N,S}v_\varepsilon(\xi_\delta)\big)+\a\mathcal{M}_{\lambda,\Lambda}^+\big(D^2_{\mathbbm{H}^N,S}\Phi(\xi_\delta)\big)\\
	&\qquad\qquad+\a\Lambda\delta\left((2N-2)\bigg(4\sum_{i=1}^N(x^i_\delta-x^i_0)^2+(y^i_\delta-y^i_0)^2\bigg)+2\bigg(12\sum_{i=1}^N(x_\delta^i-x^i_0)^2+(y_\delta^i-y^i_0)^2\bigg)\right)\\
	&\qquad\qquad-\b(-\Delta_{\mathbbm{H}^N})^su^\varepsilon(\xi_\delta),
\end{align}
where $(x_\delta^1,\dots,x_\delta^N,y_\delta^1,\dots,y_\delta^N)$ and $(x_0^1,\dots,x_0^N,y_0^1,\dots,y_0^N)$ are the first $2N$ coordinates of $\xi_\delta$ and $\xi_0$, respectively. Next, since $w(\xi\,o\,\eta)-w(\xi)>0$ for $\xi\,o\,\eta\in B_r^c(\xi)$ and $r_1$ small enough so we have by \eqref{eq 3.20}, 
\begin{align}\label{eq 3.21}
	\qquad-(-\Delta_{\mathbbm{H}^N})^su^\varepsilon(\xi_\delta)&=-(-\Delta_{\mathbbm{H}^N})^sv_\varepsilon(\xi_\delta)+(-\Delta_{\mathbbm{H}^N})^sw(\xi_\delta)-(-\Delta_{\mathbbm{H}^N})^s\Phi(\xi_\delta)-\delta(-\Delta_{\mathbbm{H}^N})^s(|\xi_0^{-1}\,o\,\xi|^4)(\xi_\delta)\\
	&\leq -(-\Delta_{\mathbbm{H}^N})^sv_\varepsilon(\xi_\delta)-\int_{\{\eta\in \mathbbm{H}^N:\,\xi_\delta\,o\,\eta\in \mathbbm{H}^N\setminus B_r^{\mathbbm{H}^N}(\xi_\delta)\}}\frac{\mathbbm{1}_{\{|\eta|_{\mathbbm{H}^N}\leq 1\}}\eta\cdot(\nabla_{\mathbbm{H}^N}w(\xi_\delta),\partial_tw(\xi_\delta))}{|\eta|^{Q+2s}_{\mathbbm{H}^N}}d\eta\\
	&\qquad-(-\Delta_{\mathbbm{H}^N})^s\Phi(\xi_\delta)-\delta(-\Delta_{\mathbbm{H}^N})^s(|\xi_0^{-1}\,o\,\xi|_{\mathbbm{H}^N}^4)(\xi_\delta)\\
	&=-(-\Delta_{\mathbbm{H}^N})^sv_\varepsilon(\xi_\delta)-\int_{\{\eta\in \mathbbm{H}^N:\,r\leq |\eta|_{\mathbbm{H}^N}\leq 1\}}\frac{\eta\cdot(\nabla_{\mathbbm{H}^N}w(\xi_\delta),\partial_tw(\xi_\delta))}{|\eta|^{Q+2s}_{\mathbbm{H}^N}}d\eta\\
	&\qquad-(-\Delta_{\mathbbm{H}^N})^s\Phi(\xi_\delta)-\delta(-\Delta_{\mathbbm{H}^N})^s\big(|\xi_0^{-1}\,o\,\xi|_{\mathbbm{H}^N}^4\big)(\xi_\delta).
\end{align}
Further, using \eqref{eq 3.21} together with \eqref{eq 3.20} yields
\begin{align}
	f(\xi_\delta)-d_\varepsilon&\leq \a\mathcal{M}_{\lambda,\Lambda}^+\big(D^2_{\mathbbm{H}^N,S}v_\varepsilon(\xi_\delta)\big)-\b(-\Delta_{\mathbbm{H}^N})^sv_\varepsilon(\xi_\delta)+\a\mathcal{M}_{\lambda,\Lambda}^+\big(D^2_{\mathbbm{H}^N,S}\Phi(\xi_\delta)\big)-\b(-\Delta_{\mathbbm{H}^N})^s\Phi(\xi_\delta)\\
	&\qquad +\a\Lambda\delta\left((2N-2)\bigg(4\sum_{i=1}^N(x_\delta^i-x^i_0)^2+(y_\delta^i-y^i_0)^2\bigg)+2\bigg(12\sum_{i=1}^N(x_\delta^i-x^i_0)^2+(y_\delta^i-y^i_0)^2\bigg)\right)\\
	&\qquad -\frac{\b}{2}c(N,s)\int_{\{\eta\in \mathbbm{H}^N:\,r\leq |\eta|_{\mathbbm{H}^N}\leq 1\}}\frac{\eta\cdot(\nabla_{\mathbbm{H}^N}w(\xi_\delta),\partial_tw(\xi_\delta))}{|\eta|^{Q+2s}_{\mathbbm{H}^N}}d\eta-\b\delta(-\Delta_{\mathbbm{H}^N})^s(|\xi_0^{-1}\,o\,\xi|_{\mathbbm{H}^N}^4)(\xi_\delta)\\
	&\leq g(\xi_\delta)+d_\varepsilon+\a\Lambda\delta\left((2N-2)\bigg(4\sum_{i=1}^N(x_\delta^i-x^i_0)^2+(y_\delta^i-y^i_0)^2\bigg)+2\bigg(12\sum_{i=1}^N(x^i_\delta-x^i_0)^2+(y^i_\delta-y^i_0)^2\bigg)\right)\\
	&\qquad-\frac{\b}{2}c(N,s)\int_{\{\eta\in \mathbbm{H}^N:\,r\leq |\eta|_{\mathbbm{H}^N}\leq 1\}}\frac{\eta\cdot(\nabla_{\mathbbm{H}^N}w(\xi_\delta),\partial_tw(\xi_\delta))}{|\eta|^{Q+2s}_{\mathbbm{H}^N}}d\eta-\delta(-\Delta_{\mathbbm{H}^N})^s(|\xi_0^{-1}\,o\,\xi|_{\mathbbm{H}^N}^4)(\xi_\delta)\\
	&\qquad+\a\mathcal{M}_{\lambda,\Lambda}^+\big(D^2_{\mathbbm{H}^N,S}\Phi(\xi_\delta)\big)-\b(-\Delta_{\mathbbm{H}^N})^s\Phi(\xi_\delta) \text{ (using \eqref{2})}.
\end{align}
It is easy to observe that $\xi_\delta\longrightarrow \xi_0$ and $\nabla\Gamma_w(\xi_\delta)\longrightarrow\nabla\Gamma_w(\xi_0)=0$ as $r_1\longrightarrow 0$. Now, letting $\delta\longrightarrow 0$ gives
\begin{align}
	f(\xi_0)-d_\varepsilon\leq g(\xi_0)+d_\varepsilon+\a\mathcal{M}_{\lambda,\Lambda}^+\big(D^2_{\mathbbm{H}^N,S}\Phi(\xi_0)\big)-\b(-\Delta_{\mathbbm{H}^N})^s\Phi(\xi_0).
\end{align}
In other words,
\begin{align}
	\a\mathcal{M}_{\lambda,\Lambda}^+\big(D^2_{\mathbbm{H}^N,S}\Phi(\xi_0)\big)-\b(-\Delta_{\mathbbm{H}^N})^s\Phi(\xi_0)\geq f(\xi_0)-g(\xi_0)-2d_\varepsilon.
\end{align}
The above equation clearly implies that
\begin{align}
	\a\mathcal{M}_{\lambda,\Lambda}^+\big(D^2_{\mathbbm{H}^N,S}(u^\varepsilon-v_\varepsilon)\big)-\b(-\Delta_{\mathbbm{H}^N})^s\Phi(u^\varepsilon-v_\varepsilon)\geq f-g-2d_\varepsilon \text{ in }\Omega_1
\end{align}
in the viscosity sense. Finally, letting $\varepsilon\longrightarrow 0$ along with using Theorem \ref{stability} yields the claim.\qed\\

In order to derive comparison principle, we state and prove the following lemma. We mention that a similar lemma has been proven in the Euclidean setting (see Lemma 5.5 \cite{Mou}).

\begin{lem}\label{Mou}
	There exists a $C^2(\overline{\Omega})\cap C_b(\mathbbm{H}^N)$ function $\varphi_h$ such that
\begin{align}
		\a\mathcal{M}_{\lambda,\Lambda}^+\big(D^2_{\mathbbm{H}^N,S}\varphi_h\big)-\b(-\Delta_{\mathbbm{H}^N})^s\varphi_h\leq -1, \text{ in }\Omega.
\end{align}
\end{lem}

\noindent \textbf{Proof of Lemma \ref{Mou}.} Let diam$(\Omega)=R.$ We may assume without of loss of generality that $\Omega\subset B^{R}_{\mathbbm{H}^N}(\xi_R),$ where $\xi_R\coloneqq(2R,0,\dots,0).$ Let us define
\begin{align}
	\varphi_h(\xi)=\begin{cases}
		2-e^{-C\xi_1} &\text{ for }\xi_1\geq 0\\
		\frac{1}{2}+\frac{1}{4}\left(\frac{1}{1-C\xi_1}\right)+\frac{1}{4}\big(\text{sin } 3C\xi_1+ \text{cos }\sqrt{6}\,C\xi_1\big) &\text{ for } \xi_1<0,
	\end{cases}
\end{align}
where $\xi_1$ is the first coordinate of $\xi=(\xi_1,\dots,\xi_N,\xi_{N+1},\dots,\xi_{2N},t)\in \Omega$ and $C>0$ is some constant.
It is easy to calculate for $\xi_1>0,$
\begin{align}
	\partial_{x_i}\varphi_h=\begin{cases}
		Ce^{-C\xi_1} &\text{ if }i=1,\\
		0 &\text{ if }2\leq i\leq 2N,
	\end{cases}
\end{align}
and $\partial_t \varphi_h=0.$ Also, for $\xi_1>0$,
\begin{align}
	\partial^2_{x_ix_j}\varphi_h=\begin{cases}
		-C^2e^{-C\xi_1} &\text{ if }i=j=1\\
		0 &\text{ otherwise, }
	\end{cases}
\end{align}
and $\partial^2_{tx_i}=0=\partial^2_{tt}$ for $i=1,2,\dots,2N.$ 
Using this, we first compute $\nabla_{\mathbbm{H}^N}\varphi_h(\xi)$ and $D^2_{\mathbbm{H}^N,S}\varphi_h(\xi)$ as follows:
\begin{align}
	\nabla_{\mathbbm{H}^N}\varphi_h(\xi)=\sigma(\xi)\nabla\varphi_h(\xi)=
	\begin{bmatrix}
		I_N &0_N &2y\\
		0_N &I_N &-2x
	\end{bmatrix}_{2N\times (2N+1)}
	\begin{bmatrix}
		Ce^{-C\xi_1}\\
		0\\
		\vdots\\
		0
	\end{bmatrix}_{(2N+1)\times 1}=
	\begin{bmatrix}
		Ce^{-C\xi_1}\\
		0\\
		\vdots\\
		0
	\end{bmatrix}_{2N\times 1}.
\end{align} 
Also,
\begin{align}
	D^2_{\mathbbm{H}^N}\varphi_h(\xi)&=\sigma(\xi)D^2\varphi_h\sigma^T(\xi)\\
	&=\begin{bmatrix}
		I_N &0_N &2y\\
		0_N &I_N &-2x
	\end{bmatrix}_{2N\times(2N+1)}
	\begin{bmatrix}
		-C^2e^{-C\xi_1} &0 &\dots &0\\
		0 &0 &\dots &0\\
		\vdots &\vdots &\ddots &\vdots \\ 
		0 &0 &\dots &0
	\end{bmatrix}_{(2N+1)\times(2N+1)}
	\begin{bmatrix}
		I_N &0_N\\
		0_N &I_N\\
		2y &-2x
	\end{bmatrix}_{(2N+1)\times 2N}\\
	&=\begin{bmatrix}
		-C^2e^{-C\xi_1} &0 &\dots &0\\
		0 &0 &\dots &0\\
		\vdots &\vdots &\ddots &\vdots \\ 
		0 &0 &\dots &0
	\end{bmatrix}_{2N\times 2N}.
\end{align}
Now, let $\delta=\displaystyle{\min\{1,R\}}.$ Then for any $\xi\in\Omega,$ we get
\begin{align}\label{111}
	\a\mathcal{M}_{\lambda,\Lambda}^+\big(D^2_{\mathbbm{H}^N,S}\varphi_h(\xi)\big)&-\b(-\Delta_{\mathbbm{H}^N})^s\varphi_h(\xi)\\
	&=-\lambda\a C^2e^{-C\xi_1}+\frac{\b}{2}c(N,s)\int_{\mathbbm{H}^N}\frac{\varphi_h(\xi\,o\,\eta)-\varphi_h(\xi)-\mathbbm{1}_{\{|\eta|_{\mathbbm{H}^N}\leq 1\}}\eta\cdot(\nabla_{\mathbbm{H}^N}\varphi_h(\xi),\partial_t\varphi_h(\xi))}{|\eta|_{\mathbbm{H}^N}^{Q+2s}}d\eta\\
	&=-\lambda \a C^2e^{-C\xi_1}+\frac{\b}{2}c(N,s)\int_{B_\delta^{\mathbbm{H}^N}}\frac{\varphi_h(\xi\,o\,\eta)-\varphi_h(\xi)-\mathbbm{1}_{\{|\eta|_{\mathbbm{H}^N}\leq 1\}}\eta\cdot(\nabla_{\mathbbm{H}^N}\varphi_h(\xi),\partial_t\varphi_h(\xi))}{|\eta|_{\mathbbm{H}^N}^{Q+2s}}d\eta\\
	&\quad+\frac{\b}{2}c(N,s)\int_{\mathcal{C}B_\delta^{\mathbbm{H}^N}\cap\{\eta:\eta_1\leq 0\}}\frac{\varphi_h(\xi\,o\,\eta)-\varphi_h(\xi)-\mathbbm{1}_{\{|\eta|_{\mathbbm{H}^N}\leq 1\}}\eta\cdot(\nabla_{\mathbbm{H}^N}\varphi_h(\xi),\partial_t\varphi_h(\xi))}{|\eta|_{\mathbbm{H}^N}^{Q+2s}}d\eta\\
	&\quad+\frac{\b}{2}c(N,s)\int_{B^{\mathbbm{H}^N}_1\cap \,\mathcal{C}B_\delta^{\mathbbm{H}^N}\cap\,\{\eta:\eta_1>0\}}\frac{\varphi_h(\xi\,o\,\eta)-\varphi_h(\xi)-\mathbbm{1}_{\{|\eta|_{\mathbbm{H}^N}\leq 1\}}\eta\cdot(\nabla_{\mathbbm{H}^N}\varphi_h(\xi),\partial_t\varphi_h(\xi))}{|\eta|_{\mathbbm{H}^N}^{Q+2s}}d\eta\\
	&\quad+\frac{\b}{2}c(N,s)\int_{\mathcal{C}B^{\mathbbm{H}^N}_1\cap\,\{\eta:\eta_1>0\}}\frac{\varphi_h(\xi\,o\,\eta)-\varphi_h(\xi)-\mathbbm{1}_{\{|\eta|_{\mathbbm{H}^N}\leq 1\}}\eta\cdot(\nabla_{\mathbbm{H}^N}\varphi_h(\xi),\partial_t\varphi_h(\xi))}{|\eta|_{\mathbbm{H}^N}^{Q+2s}}d\eta,
\end{align}
where $\mathcal{C}B_r^{\mathbbm{H}^N}$ denotes the complement of $B_r^{\mathbbm{H}^N}$ in $\mathbbm{H}^N,$ i.e., $\mathbbm{H}^N\setminus B_r^{\mathbbm{H}^N}.$ We can re-write the above equation as follows:
\begin{align}
	\a\mathcal{M}_{\lambda,\Lambda}^+&\big(D^2_{\mathbbm{H}^N,S}\varphi_h(\xi)\big)-\b(-\Delta_{\mathbbm{H}^N})^s\varphi_h(\xi)\\
	&=-\lambda \a C^2e^{-C\xi_1}+\frac{\b}{2}c(N,s)\int_{B_\delta^{\mathbbm{H}^N}}\frac{\varphi_h(\xi\,o\,\eta)-\varphi_h(\xi)-\eta\cdot(\nabla_{\mathbbm{H}^N}\varphi_h(\xi),\partial_t\varphi_h(\xi))}{|\eta|_{\mathbbm{H}^N}^{Q+2s}}d\eta\\
	&\qquad+\frac{\b}{2}c(N,s)\int_{\mathcal{C}B_\delta^{\mathbbm{H}^N}\cap\{\eta:\eta_1\leq 0\}}\frac{\varphi_h(\xi\,o\,\eta)-\varphi_h(\xi)}{|\eta|_{\mathbbm{H}^N}^{Q+2s}}d\eta+\frac{\b}{2}c(N,s)\int_{B_1^{\mathbbm{H}^N}\cap \,\mathcal{C}B_\delta^{\mathbbm{H}^N}\cap\,\{\eta:\eta_1>0\}}\frac{\varphi_h(\xi\,o\,\eta)-\varphi_h(\xi)}{|\eta|_{\mathbbm{H}^N}^{Q+2s}}d\eta\\
	&\qquad+\frac{\b}{2}c(N,s)\int_{B^{\mathbbm{H}^N}_1\cap\mathcal{C}B_\delta^{\mathbbm{H}^N}}-\frac{\eta\cdot(\nabla_{\mathbbm{H}^N}\varphi_h(\xi),\partial_t\varphi_h(\xi))}{|\eta|_{\mathbbm{H}^N}^{Q+2s}}d\eta+\frac{\b}{2}c(N,s)\int_{\mathcal{C}B^{\mathbbm{H}^N}_1\cap\,\{\eta:\eta_1>0\}}\frac{\varphi_h(\xi\,o\,\eta)-\varphi_h(\xi)}{|\eta|_{\mathbbm{H}^N}^{Q+2s}}d\eta.
\end{align}
It further gives
\begin{align}\label{112}
	\qquad\quad
	\a\mathcal{M}_{\lambda,\Lambda}^+(D^2_{\mathbbm{H}^N,S}\varphi_h(\xi))-\b(-\Delta_{\mathbbm{H}^N})^s\varphi_h(\xi)&=-\lambda \a C^2e^{-C\xi_1}\\
	&\qquad+\frac{\b}{2}c(N,s)\int_{B_\delta^{\mathbbm{H}^N}}\frac{\big(2-e^{-C(\xi_1+\eta_1)}\big)-\big(2-e^{-C\xi_1}\big)-\eta_1Ce^{-C\xi_1}}{|\eta|_{\mathbbm{H}^N}^{Q+2s}}d\eta\\
	&\qquad+\frac{\b}{2}c(N,s)\int_{\mathcal{C}B_\delta^{\mathbbm{H}^N}\cap\{\eta:\eta_1\leq 0\}}\frac{\big(2-e^{-C(\xi_1+\eta_1)}\big)-\big(2-e^{-C\xi_1}\big)}{|\eta|_{\mathbbm{H}^N}^{Q+2s}}d\eta\\
	&\qquad+\frac{\b}{2}c(N,s)\int_{B_1^{\mathbbm{H}^N}\cap \,\mathcal{C}B_\delta^{\mathbbm{H}^N}\cap\,\{\eta:\eta_1>0\}}\frac{\big(2-e^{-C(\xi_1+\eta_1)}\big)-\big(2-e^{-C\xi_1}\big)}{|\eta|_{\mathbbm{H}^N}^{Q+2s}}d\eta\\
	&\qquad+\frac{\b}{2}c(N,s)\int_{B^{\mathbbm{H}^N}_1\cap\mathcal{C}B_\delta^{\mathbbm{H}^N}}-\frac{\eta_1Ce^{-C\xi_1}}{|\eta|_{\mathbbm{H}^N}^{Q+2s}}d\eta\\
	&\qquad+\frac{\b}{2}c(N,s)\int_{\mathcal{C}B_1^{\mathbbm{H}^N}\cap\,\{\eta:\eta_1>0\}}\frac{\big(2-e^{-C(\xi_1+\eta_1)}\big)-\big(2-e^{-C\xi_1}\big)}{|\eta|_{\mathbbm{H}^N}^{Q+2s}}d\eta.
\end{align}
Similarly,
\begin{align}\label{112}
\qquad\quad
\a\mathcal{M}_{\lambda,\Lambda}^+(D^2_{\mathbbm{H}^N,S}\varphi_h(\xi))&-\b(-\Delta_{\mathbbm{H}^N})^s\varphi_h(\xi)\\
&=-\lambda\alpha C^2e^{-C\xi_1}+\frac{\b}{2}c(N,s)\int_{B_\delta^{\mathbbm{H}^N}}\frac{\big(-e^{-C(\xi_1+\eta_1)}+e^{-C\xi_1}\big)-\eta_1Ce^{-C\xi_1}}{|\eta|_{\mathbbm{H}^N}^{Q+2s}}d\eta\\
	&\qquad+\frac{\b}{2}c(N,s)\int_{\mathcal{C}B_\delta^{\mathbbm{H}^N}\cap\{\eta:\eta_1\leq 0\}}\frac{\big(-e^{-C(\xi_1+\eta_1)}+e^{-C\xi_1}\big)}{|\eta|_{\mathbbm{H}^N}^{Q+2s}}d\eta\\
	&\qquad+\frac{\b}{2}c(N,s)\int_{B_1^{\mathbbm{H}^N}\cap \,\mathcal{C}B_\delta^{\mathbbm{H}^N}\cap\,\{\eta:\eta_1>0\}}\frac{\big(-e^{-C(\xi_1+\eta_1)}+e^{-C\xi_1}\big)}{|\eta|_{\mathbbm{H}^N}^{Q+2s}}d\eta\\
	&\qquad+\frac{\b}{2}c(N,s)\int_{B^{\mathbbm{H}^N}_1\cap\mathcal{C}B_\delta^{\mathbbm{H}^N}}-\frac{\eta_1Ce^{-C\xi_1}}{|\eta|_{\mathbbm{H}^N}^{Q+2s}}d\eta+\b\int_{\mathcal{C}B_1^{\mathbbm{H}^N}\cap\,\{\eta:\eta_1>0\}}\frac{\big(-e^{-C(\xi_1+\eta_1)}+e^{-C\xi_1}\big)}{|\eta|_{\mathbbm{H}^N}^{Q+2s}}d\eta.
\end{align} 
Now, since $e^{-C\xi_1}$ is a convex function so we have
\begin{align}\label{int 1}
	0&\leq e^{-C(\xi_1+\eta_1)}-e^{-C\xi_1}+\eta\cdot(\nabla_{\mathbbm{H}^N}\varphi_h(\xi),\partial_t\varphi_h(\xi))\\
	&= e^{-C(\xi_1+\eta_1)}-e^{-C\xi_1}+C\eta_1e^{-C\xi_1}.
\end{align}
Also, if $\eta_1\leq 0$ then 
\begin{align}\label{int 2}
	e^{-C(\xi_1+\eta_1)}-e^{-C\xi_1}\geq 0.
\end{align}
Using \eqref{int 1} and \eqref{int 2} in the first and second integrals in the R.H.S. of \eqref{112}, respectively, it confers that these integrals are non-positive. We further have that
\begin{align}
	|e^{-C(\xi_1+\eta_1)}-e^{-C\xi_1}|= e^{-C\xi_1}|e^{-C\eta_1}-1|\leq Ce^{-C\xi_1}|\eta|_{\mathbbm{H}^N}.
\end{align}
It infers
\begin{align}\label{int 3}
	\left|\b\int_{B_1^{\mathbbm{H}^N}\cap \,\mathcal{C}B_\delta^{\mathbbm{H}^N}\cap\,\{\eta:\eta_1>0\}}\frac{\left(-e^{-C(\xi_1+\eta_1)}+e^{-C\xi_1}\right)}{|\eta|_{\mathbbm{H}^N}^{Q+2s}}d\eta\right|&\leq \b\int_{B_1^{\mathbbm{H}^N}\cap \,\mathcal{C}B_\delta^{\mathbbm{H}^N}}\frac{Ce^{-C\xi_1}|\eta|_{\mathbbm{H}^N}}{|\eta|_{\mathbbm{H}^N}^{Q+2s}}d\eta\\
	&\leq \b\int_{B_1^{\mathbbm{H}^N}\cap \,\mathcal{C}B_\delta^{\mathbbm{H}^N}}\frac{Ce^{-C\xi_1}|\eta|^2_{\mathbbm{H}^N}}{\delta|\eta|_{\mathbbm{H}^N}^{Q+2s}}d\eta\\
	&\leq \b e^{-C\xi_1}\frac{C}{\delta}\int_{B^{\mathbbm{H}^N}_1\cap \,\mathcal{C}B_\delta^{\mathbbm{H}^N}}\frac{|\eta|^2_{\mathbbm{H}^N}}{|\eta|_{\mathbbm{H}^N}^{Q+2s}}d\eta\\
	&\leq \b e^{-C\xi_1}\frac{C}{\delta}\int_{B_1^{\mathbbm{H}^N}}\frac{|\eta|^2_{\mathbbm{H}^N}}{|\eta|_{\mathbbm{H}^N}^{Q+2s}}d\eta.
\end{align}
It is easy to see that for $\eta_1>0,$ we have $\big|e^{-C(\xi_1+\eta_1)}-e^{-C\xi_1}\big|\leq e^{-C\xi_1}$ so
\begin{align}\label{int 5}
	\left|\b\int_{\mathcal{C}B_1^{\mathbbm{H}^N}\cap\,\{\eta:\eta_1>0\}}\frac{\big(-e^{-C(\xi_1+\eta_1)}+e^{-C\xi_1}\big)}{|\eta|_{\mathbbm{H}^N}^{Q+2s}}d\eta\right|\leq \b e^{-C\xi_1}\int_{\mathcal{C}B_1^{\mathbbm{H}^N}}\frac{1}{|\eta|_{\mathbbm{H}^N}^{Q+2s}}d\eta.
\end{align}
Also, it is easy to observe that
\begin{align}\label{int 4}
	\left|\b\int_{B^{\mathbbm{H}^N}_1\cap\mathcal{C}B_\delta^{\mathbbm{H}^N}}-\frac{\eta_1Ce^{-C\xi_1}}{|\eta|_{\mathbbm{H}^N}^{Q+2s}}d\eta\right|&\leq \b\int_{B^{\mathbbm{H}^N}_1\cap\mathcal{C}B_\delta^{\mathbbm{H}^N}}\frac{|\eta|_{\mathbbm{H}^N}Ce^{-C\xi_1}}{|\eta|_{\mathbbm{H}^N}^{Q+2s}}d\eta\\
	&\leq \b e^{-C\xi_1}\frac{C}{\delta}\int_{B^{\mathbbm{H}^N}_1}\frac{|\eta|^2_{\mathbbm{H}^N}}{|\eta|_{\mathbbm{H}^N}^{Q+2s}}d\eta.
\end{align} 
Using \eqref{int 3}, \eqref{int 5} \& \eqref{int 4} in \eqref{112}
yields that for sufficiently large enough $C>0$ we can make the L.H.S. of \eqref{112} less than $-1$, i.e.,
\begin{align}
	\a\mathcal{M}_{\lambda,\Lambda}^+\big(D^2_{\mathbbm{H}^N,S}\varphi_h\big)-\b(-\Delta_{\mathbbm{H}^N})^s\varphi_h&=\leq -1.
\end{align}
This competes the proof.\qed\\

Further, using the above lemma, we prove the comparison principle. Lemma \ref{Lemma 1.6} together with a standard approximation argument produces the following lemma.

\begin{lem}\label{comp} 
	Let $u$ be bounded in $\mathbbm{H}^N$ and USC in $\overline{\Omega}$ such that
	\begin{align}
		\a\mathcal{M}^+_{\lambda,\Lambda}\big(D^2_{\mathbbm{H}^N,S}u\big)-\b(-\Delta_{\mathbbm{H}^N})^su\geq 0 \text{ in }\Omega \text{ in the viscosity sense.}
	\end{align}
	Then
	\begin{align}
		\displaystyle{\sup_\Omega}\,u\leq \displaystyle{\sup_{\mathbbm{H}^N\setminus\Omega}u}.
	\end{align}
\end{lem}
\noindent \textbf{Proof of Lemma \ref{comp}.} By Lemma \ref{Mou}, we have a function $\varphi_h\in C^2(\overline{\Omega})\cap C_b(\mathbbm{H}^N)$ such that 
\begin{align}
	\a\mathcal{M}_{\lambda,\Lambda}^+\big(D^2_{\mathbbm{H}^N,S}\varphi_h\big)-\b(-\Delta_{\mathbbm{H}^N})^s\varphi_h\leq -1.
\end{align}
Now, let us define a function
\begin{align}
	\varphi_M(x)=M+{\mathlarger\varepsilon}\varphi_h(x)\text{ for }\mathlarger\varepsilon>0.
\end{align}
This provides
\begin{align}\label{contra}
	\a\mathcal{M}_{\lambda,\Lambda}^+\big(D^2_{\mathbbm{H}^N,S}\varphi_M\big)-\b(-\Delta_{\mathbbm{H}^N})^s\varphi_M&=\alpha\mathlarger\varepsilon\mathcal{M}_{\lambda,\Lambda}^+\big(D^2_{\mathbbm{H}^N,S}\varphi_h\big)-\b\mathlarger\varepsilon(-\Delta_{\mathbbm{H}^N})^s\varphi_h\\
	&\leq -\varepsilon \text{ in }\Omega.
\end{align}
Further, let $M_0$ be the smallest value of $M$ such that $\varphi_M\geq u.$ We aim to prove $M_0\leq \displaystyle{\sup_{{\mathbbm{H}^N}\setminus\Omega}}u$ by the method of contradiction. Let us assume that
\begin{align}
	M_0>\displaystyle{\sup_{{\mathbbm{H}^N}\setminus\Omega}}u,
\end{align}
then we have that there exists a point $\xi_0\in \Omega$ such that $u(\xi_0)=\varphi_{M_0}(\xi_0).$ It immediately implies that $\varphi_{M_0}$ touches $u$ at $\xi_0$ from above and by the definition, it infers 
\begin{align}
	\a\mathcal{M}_{\lambda,\Lambda}^+\big(D^2_{\mathbbm{H}^N,S}\varphi_{M_0}\big)(\xi_0)-\b(-\Delta_{\mathbbm{H}^N})^s\varphi_{M_0}(\xi_0)\geq 0.
\end{align}
This contradicts \eqref{contra}. Thus, we have
\begin{align}
	M_0\leq \displaystyle{\sup_{{\mathbbm{H}^N}\setminus\Omega}}u,
\end{align}
 which further entails that
 \begin{align}
 	u(\xi)&\leq \varphi_{M_0}(\xi_0)\\
 	 &\leq M_0+\mathlarger\varepsilon\,\displaystyle{\sup_{\mathbbm{H}^N}\varphi_h}\\
 	 &\leq \displaystyle{\sup_{\mathbbm{H}^N\setminus\Omega}u}+\mathlarger\varepsilon\,\displaystyle{\sup_{\mathbbm{H}^N}}\,\varphi_h,\, \xi\in \mathbbm{H}^N.
 \end{align}
Finally, letting $\mathlarger\varepsilon\longrightarrow 0$ offers the claim.\qed\\

\noindent \textbf{Proof of Theorem \ref{Holder}.} We follow the similar arguments as in the Euclidean setting, see for instance, \cite{Silvestre, Lara, Silvestre 2}. Without loss of generality, we may assume that $\|u\|_{\infty,\mathbbm{H}^N}=1.$ We first show the existence of a universal constant $0<\delta<1$ such that
\begin{align}
	\underset{{B_{2^{-k}}^{\mathbbm{H}^N}}}{\text{osc}}u\leq 2(1-\delta)^k.
\end{align}
We use the principle of mathematical induction to show the claim. Now, since $\|u\|_{\infty,\mathbbm{H}^N}=1$ so the claim holds trivially for $k\leq 0.$ Let the above inequality holds true up to some $k$, we need to show that it also holds for $k+1.$ For this, consider a function
\begin{align}
	v(\xi)\coloneqq\frac{u(2^{-k}\xi)}{(1-\delta)^k
	}-a_k\,,
\end{align}
where $a_k$ is a constant chosen such that $-\frac{1}{2}\leq v\leq \frac{1}{2}$ in $B_1^{\mathbbm{H}^N}.$ Also, by the induction hypothesis, we have
\begin{align}
	\underset{{B_{2^{-j}}^{\mathbbm{H}^N}}}{\text{osc}}u\leq 2(1-\delta)^j \text{ for all }j\leq k,
\end{align} 
i.e.,
\begin{align}
	\underset{{B_{2^{j}}^{\mathbbm{H}^N}}}{\text{osc}}u\leq 2(1-\delta)^{-j} \text{ for all }j\geq 0.
\end{align}
We aim to show that
\begin{align}
	\underset{{B_{\frac{1}{2}}^{\mathbbm{H}^N}}}{\text{osc}}v\leq \underset{{B_{1}^{\mathbbm{H}^N}}}{\text{osc}}v.
\end{align}
One may show this either by showing that supremum of $v$ in $B_{\frac{1}{2}}^{\mathbbm{H}^N}$ is smaller than that in $B_{1}^{\mathbbm{H}^N}$ or the infimum is larger. It is trivial that either $v\geq 0$ or $v\leq 0$ for atleast half of the points (in measure) in $B_1^{\mathbbm{H}^N}$. Without loss of generality, we may assume that 
\begin{align}
	|N|\geq \frac{1}{2}\big|B_1^{\mathbbm{H}^N}\big|,
\end{align}   
where $N\coloneqq\{v\leq 0\}\cap B_1^{\mathbbm{H}^N}.$ We show now that the induction hypothesis implies that  
\begin{align}
	v(\xi)\leq \big(2|\xi|_{{\mathbbm{H}^N}}\big)^\a-\frac{1}{2} \text{ for all } \xi\in B_1^{\mathbbm{H}^N}.
\end{align}
Let $\xi\in {B_{2^{j+1}}^{\mathbbm{H}^N}}\setminus{B_{2^{j}}^{\mathbbm{H}^N}}$ for some $j>0.$ That is 
\begin{align}
	2^{j+1}>|\xi|_{\mathbbm{H}^N}>2^j.
\end{align}
In other words,
\begin{align}
	2^{-k+j+1}>|2^{-k}\xi|_{\mathbbm{H}^N}>2^{-k+j}.
\end{align}
\begin{align}
	v(x)=(1-\delta)^{-k}u(2^{-k}x)-a_k.
\end{align}
For $\xi\in {B_{1}^{\mathbbm{H}^N}},$
\begin{align}
	-\frac{1}{2}\leq (1-\delta)^{-k}u(2^{-k}\xi)-a_k\leq \frac{1}{2}.
\end{align}
Now, for $\xi\notin {B_{1}^{\mathbbm{H}^N}}$, we have for any index $j\geq 0$,
\begin{align}\label{eq 3.23}
	v(\xi)\leq \big(2|\xi|_{\mathbbm{H}^N}\big)^{\c}-\frac{1}{2},
\end{align}
where $\c$ is a number such that $(1-\delta)=2^{-\c}.$
Also, for any $\xi\in B_1^{\mathbbm{H}^N},$
\begin{align}
	\a\mathcal{M}_{\lambda,\Lambda}^+\big(D^2_{\mathbbm{H}^N}v\big)-\b(-\Delta_{\mathbbm{H}^N})^sv=(1-\delta)^k\left(\a\mathcal{M}_{\lambda,\Lambda}^+\big(D^2_{\mathbbm{H}^N}u\big)-\b(-\Delta_{\mathbbm{H}^N})^su\right)=0.	
\end{align}
Now, we show that the following three points:
\begin{enumerate}
\item[(i)]$v(\xi)\leq \big(2|\xi|_{\mathbbm{H}^N}\big)^\c-\frac{1}{2}$ for $\xi\notin B_1^{\mathbbm{H}^N}$
\item[(ii)]$|N|=\big|\{v\leq 0\}\cap B_1^{\mathbbm{H}^N}\big|\geq \frac{1}{2}\big|B_1^{\mathbbm{H}^N}\big|$
\item[(iii)] $\a\mathcal{M}_{\lambda,\Lambda}^+(D^2_{\mathbbm{H}^N}v)-\b(-\Delta_{\mathbbm{H}^N})^sv=0$ in $B_1^{\mathbbm{H}^N}$
\end{enumerate}
imply that $v\leq \big({1}/{2}-\delta\big)$ in $B_{{1}/{2}}^{\mathbbm{H}^N}.$ Let us assume the contrary. Let $v(\xi)>\left({1}/{2}-\delta\right)$ for some $\xi\in B_{{1}/{2}}^{\mathbbm{H}^N}.$ Consider a smooth radial function whose support is contained in $B_{{3}/{4}}^{\mathbbm{H}^N}$ and $\rho\equiv 1$ in $B_{{1}/{2}}^{\mathbbm{H}^N}.$ It immediately gives that $v+\delta\rho$ attains a local maximum at some point $\xi_0\in B_{{3}/{4}}^{\mathbbm{H}^N}$ such that $\big(v+\delta\rho
\big)(\xi_0)>{1}/{2},$\, i.e.,\, $\displaystyle{\max_{B_1^{\mathbbm{H}^N}}}(v+\delta\rho)=(v+\delta\rho)(\xi_0)>{1}/{2}.$
Let us evaluate $\mathcal{L}(v+\delta\rho)$, which further entails getting a contradiction. 
\begin{align}
	\mathcal{L}\big(v+\delta\rho\big)(\xi_0)&=\a\mathcal{M}_{\lambda,\Lambda}^+\big(D^2_{\mathbbm{H}^N}(v+\delta\rho)\big)(\xi_0)-\b(-\Delta_{\mathbbm{H}^N})^s(v+\delta\rho)(\xi_0)\\
	&\geq \a\mathcal{M}_{\lambda,\Lambda}^+\big(D^2_{\mathbbm{H}^N}v\big)(\xi_0)+\delta\a\mathcal{M}_{\lambda,\Lambda}^-\big(D^2_{\mathbbm{H}^N}\rho\big)(\xi_0)-\b(-\Delta_{\mathbbm{H}^N})^sv(\xi_0)-\delta\b(-\Delta_{\mathbbm{H}^N})^s\rho(\xi_0)\\
	&=\a\mathcal{M}_{\lambda,\Lambda}^+\big(D^2_{\mathbbm{H}^N}v\big)(\xi_0)-\b(-\Delta_{\mathbbm{H}^N})^sv(\xi_0)+\delta\a\mathcal{M}_{\lambda,\Lambda}^+\big(D^2_{\mathbbm{H}^N}\rho\big)(\xi_0)-\delta\b(-\Delta_{\mathbbm{H}^N})^s\rho(\xi_0)\\
	&=\delta\a\mathcal{M}_{\lambda,\Lambda}^+\big(D^2_{\mathbbm{H}^N}\rho\big)(\xi_0)-\delta\b(-\Delta_{\mathbbm{H}^N})^s\rho(\xi_0)\\
	&\geq \delta\,\displaystyle{\min_{\xi\in B^{\mathbbm{H}^N}_\frac{3}{4}}}\big(\a\mathcal{M}_{\lambda,\Lambda}^+\big(D^2_{\mathbbm{H}^N}\rho\big)(\xi)-\b(-\Delta_{\mathbbm{H}^N})^s\rho(\xi)\big).
\end{align}
On the other hand,
\begin{align}\label{eq 3.24}
	\\\mathcal{L}&(v+\delta\rho)(\xi_0)\\
	&=\a\mathcal{M}_{\lambda,\Lambda}^+(D^2_{\mathbbm{H}^N}(v+\delta\rho))(\xi_0)-\b(-\Delta_{\mathbbm{H}^N})^s(v+\delta\rho)(\xi_0)\\
	&\leq \a\mathcal{M}_{\lambda,\Lambda}^+(D^2_{\mathbbm{H}^N}(v+\delta\rho))(\xi_0)+\frac{\b}{2}c(N,s)\int_{\mathbbm{H}^N}\frac{(v+\delta\rho)(\xi_0\,o\,\eta)+(v+\delta\rho)(\xi_0\,o\,\eta^{-1})-2(v+\delta\rho)(\xi_0)}{|\eta|^{Q+2s}_{\mathbbm{H}^N}}d\eta\\
	&=\a\mathcal{M}_{\lambda,\Lambda}^+(D^2_{\mathbbm{H}^N}(v+\delta\rho))(\xi_0)\\
	&\qquad\qquad+
	\frac{\b}{2}c(N,s)\int_{\mathbbm{H}^N}\frac{(v+\delta\rho)(\xi_0\,o\,\eta)-(v+\delta\rho)(\xi_0)-\mathbbm{1}_{\{|\eta|_{\mathbbm{H}^N}\leq 1\}}\eta.\big(\nabla_{\mathbbm{H}^N}(v+\delta\rho)(\xi_0),\,\partial_t(v+\delta\rho)(\xi_0)\big)}{|\eta|^{Q+2s}_{\mathbbm{H}^N}}d\eta\\
	&=\a\mathcal{M}_{\lambda,\Lambda}^+(D^2_{\mathbbm{H}^N}(v+\delta\rho))(\xi_0)\\
	&\,\,+
	\frac{\b}{2}c(N,s)\int_{\big\{\eta\in\mathbbm{H}^N;\xi_0\,o\,\eta\in B_1^{\mathbbm{H}^N}\big\}}\frac{(v+\delta\rho)(\xi_0\,o\,\eta)-(v+\delta\rho)(\xi_0)-\mathbbm{1}_{\{|\eta|_{\mathbbm{H}^N}\leq 1\}}\eta.\big(\nabla_{\mathbbm{H}^N}(v+\delta\rho)(\xi_0),\,\partial_t(v+\delta\rho)(\xi_0)\big)}{|\eta|^{Q+2s}_{\mathbbm{H}^N}}d\eta\\
	&\,\,+
	\frac{\b}{2}c(N,s)\int_{\big\{\eta\in\mathbbm{H}^N:\xi_0\,o\,\eta\notin B_1^{\mathbbm{H}^N}\big\}}\frac{(v+\delta\rho)(\xi_0\,o\,\eta)-(v+\delta\rho)(\xi_0)-\mathbbm{1}_{\{|\eta|_{\mathbbm{H}^N}\leq 1\}}\eta.\big(\nabla_{\mathbbm{H}^N}(v+\delta\rho)(\xi_0),\,\partial_t(v+\delta\rho)(\xi_0)\big)}{|\eta|^{Q+2s}_{\mathbbm{H}^N}}d\eta.
\end{align}
Further, using the fact that $v+\delta\rho$ has a local maximum at $\xi_0$ yields
\begin{align}
	\mathcal{L}(v+\delta\rho)(\xi_0)&=\a\mathcal{M}_{\lambda,\Lambda}^+(D^2_{\mathbbm{H}^N}(v+\delta\rho))(\xi_0)\\
	&\quad\qquad+\frac{\b}{2}c(N,s)\int_{\big\{\eta\in\mathbbm{H}^N:\xi_0\,o\,\eta\in B_1^{\mathbbm{H}^N}\big\}}\frac{(v+\delta\rho)(\xi_0\,o\,\eta)-(v+\delta\rho)(\xi_0)}{|\eta|^{Q+2s}_{\mathbbm{H}^N}}d\eta\\
	&\qquad\quad+
	\frac{\b}{2}c(N,s)\int_{\big\{\eta\in\mathbbm{H}^N:\xi_0\,o\,\eta\notin B_1^{\mathbbm{H}^N}\big\}}\frac{(v+\delta\rho)(\xi_0\,o\,\eta)-(v+\delta\rho)(\xi_0)}{|\eta|^{Q+2s}_{\mathbbm{H}^N}}d\eta\\
	&\leq \frac{\b}{2}c(N,s)\int_{\big\{\eta\in\mathbbm{H}^N:\xi_0\,o\,\eta\in B_1^{\mathbbm{H}^N}\big\}}\frac{(v+\delta\rho)(\xi_0\,o\,\eta)-(v+\delta\rho)(\xi_0)}{|\eta|^{Q+2s}_{\mathbbm{H}^N}}d\eta\\
	&\qquad\quad+\frac{\b}{2}c(N,s)\int_{\big\{\eta\in\mathbbm{H}^N:\xi_0\,o\,\eta\notin B_1^{\mathbbm{H}^N}\big\}}\frac{(v+\delta\rho)(\xi_0\,o\,\eta)-(v+\delta\rho)(\xi_0)}{|\eta|^{Q+2s}_{\mathbbm{H}^N}}d\eta.
\end{align}
This also grants the non-positivity of integrand in the last integral. Moreover, we have 
\begin{align}
	\frac{\b}{2}c(N,s)\int_{\big\{\eta\in\mathbbm{H}^N:\xi_0\,o\,\eta\in B_1^{\mathbbm{H}^N}\big\}}&\frac{(v+\delta\rho)(\xi_0\,o\,\eta)-(v+\delta\rho)(\xi_0)}{|\eta|^{Q+2s}_{\mathbbm{H}^N}}d\eta\\
	&=\frac{\b}{2}c(N,s)\int_{\big\{\eta\in\mathbbm{H}^N:\xi_0\,o\,\eta\in N\big\}}\frac{(v+\delta\rho)(\xi_0\,o\,\eta)-(v+\delta\rho)(\xi)}{|\eta|^{Q+2s}_{\mathbbm{H}^N}}d\eta\\
	&\qquad+\frac{\b}{2}c(N,s)\int_{\big\{\eta\in\mathbbm{H}^N:\xi_0\,o\,\eta\in B_1^{\mathbbm{H}^N}\setminus N\big\}}\frac{(v+\delta\rho)(\xi_0\,o\,\eta)-(v+\delta\rho)(\xi_0)}{|\eta|^{Q+2s}_{\mathbbm{H}^N}}d\eta\\
	&\leq \frac{\b}{2}c(N,s)\int_{\big\{\eta\in\mathbbm{H}^N:\xi_0\,o\,\eta\in N\big\}}\frac{(v+\delta\rho)(\xi_0\,o\,\eta)-(v+\delta\rho)(\xi_0)}{|\eta|^{Q+2s}_{\mathbbm{H}^N}}d\eta.
\end{align}
This further infers
	\begin{align}		
\frac{\b}{2}c(N,s)\int_{\big\{\eta\in\mathbbm{H}^N:\xi_0\,o\,\eta\in B_1^{\mathbbm{H}^N}\big\}}&\frac{(v+\delta\rho)(\xi_0\,o\,\eta)-(v+\delta\rho)(\xi_0)}{|\eta|^{Q+2s}_{\mathbbm{H}^N}}d\eta\\	
&\leq \frac{\b}{2}c(N,s)\int_{\big\{\eta\in\mathbbm{H}^N:\xi_0\,o\,\eta\in N\big\}}\frac{\delta\rho(\xi_0\,o\,\eta)-(v+\delta\rho)(\xi_0)}{|\eta|^{Q+2s}_{\mathbbm{H}^N}}d\eta\\
	&<\frac{\b}{2}c(N,s)\int_{\big\{\eta\in\mathbbm{H}^N:\xi_0\,o\,\eta\in N\big\}}\frac{\big(\delta\rho(\xi_0\,o\,\eta)-\frac{1}{2}\big)}{|\eta|^{Q+2s}_{\mathbbm{H}^N}}d\eta\\
	&\leq \frac{\b}{2}c(N,s)\int_{\big\{\eta\in\mathbbm{H}^N:\xi_0\,o\,\eta\in N\big\}}\frac{\big(\delta M-\frac{1}{2}\big)}{|\eta|^{Q+2s}_{\mathbbm{H}^N}}d\eta,
\end{align}
for $M=\displaystyle{\max_{B^{\mathbbm{H}^N}_{{3}/{4}}}}\,\rho.$ Taking $\delta<{1}/{2M}$ along with using $|N|\geq {\big|B_1^{\mathbbm{H}^N}\big|}/{2},$ we get
\begin{align}\label{eq 3.25}
	\frac{\b}{2}c(N,s)\int_{\big\{\eta\in\mathbbm{H}^N:\xi_0\,o\,\eta\in N\big\}}\frac{(v+\delta\rho)(\xi_0\,o\,\eta)-(v+\delta\rho)(\xi_0)}{|\eta|^{Q+2s}_{\mathbbm{H}^N}}d\eta\leq -C \text{ for some }C>0. 
\end{align}
Next, using \eqref{eq 3.23}, we get a bound on the integrand of first integral in the last line of \eqref{eq 3.24} . In particular, we have 
\begin{align}\label{eq 3.26}
\frac{\beta}{2}c(N,s)\int_{\big\{\eta\in\mathbbm{H}^N:\xi_0\,o\,\eta\notin B_1^{\mathbbm{H}^N}\big\}}&\frac{(v+\delta\rho)(\xi_0\,o\,\eta)-(v+\delta\rho)(\xi_0)}{|\eta|^{Q+2s}_{\mathbbm{H}^N}}d\eta\\
&=\frac{\b}{2}c(N,s)\int_{\big\{\eta\in\mathbbm{H}^N:\xi_0\,o\,\eta\notin B_1^{\mathbbm{H}^N}\big\}}\frac{v(\xi_0\,o\,\eta)-(v+\delta\rho)(\xi_0)}{|\eta|^{Q+2s}_{\mathbbm{H}^N}}d\eta\\
&\leq \frac{\b}{2}c(N,s)\int_{\big\{\eta\in\mathbbm{H}^N:\xi_0\,o\,\eta\notin B_1^{\mathbbm{H}^N}\big\}}\frac{|\xi_0\,o\,\eta|^\c_{\mathbbm{H}^N}-\frac{1}{2}-\frac{1}{2}}{|\eta|^{Q+2s}_{\mathbbm{H}^N}}d\eta\\
&=\frac{\b}{2}c(N,s)\int_{\big\{\eta\in\mathbbm{H}^N:\xi_0\,o\,\eta\notin B_1^{\mathbbm{H}^N}\big\}}\frac{|\xi_0\,o\,\eta|^\c_{\mathbbm{H}^N}-1}{|\eta|^{Q+2s}_{\mathbbm{H}^N}}d\eta.
\end{align}
Further, taking small enough $\gamma$, we can make the above integral in the R.H.S. much less than $1.$ Therefore, by using together \eqref{eq 3.25} and \eqref{eq 3.26} for small enough $\delta$ and $\gamma$, we can make the L.H.S. of \eqref{eq 3.24} arbitrarily small. Hence, we can make it smaller than $$\delta \displaystyle{\min_{\xi\in B_{\frac{3}{4}}^{\mathbbm{H}^N}}\big(\a\mathcal{M}_{\lambda,\Lambda}^-(D^2_{\mathbbm{H}^N}\rho)(\xi)-\b(-\Delta_{\mathbbm{H}^N})^s\rho(\xi)\big)},$$ which yield a contradiction. Thus, we have
\begin{align}
	v\leq \frac{1}{2}-\delta \text{ in }B_{\frac{1}{2}}^{\mathbbm{H}^N}.
\end{align}
It gives
\begin{align}
	\underset{{B_{\frac{1}{2}}^{\mathbbm{H}^N}}}{\text{osc}}v\leq (1-\delta)\underset{{B_{1}^{\mathbbm{H}^N}}}{\text{osc}}v.
\end{align}
More precisely, we have
\begin{align}
	\frac{1}{(1-\delta)^k}\displaystyle{\underset{{B_{\frac{1}{2}}^{\mathbbm{H}^N}}}{\text{osc}}u(2^{-k}\xi)}=\frac{1}{(1-\delta)^k}\displaystyle{\underset{{B_{2^{-k-1}}^{\mathbbm{H}^N}}}{\text{osc}}u(\xi)}\leq \frac{(1-\delta)}{(1-\delta)^k}\displaystyle{\underset{{B_{1}^{\mathbbm{H}^N}}}{\text{osc}}u(2^{-k}\xi)}=\frac{(1-\delta)}{(1-\delta)^k}\displaystyle{\underset{{B_{2^{-k}}^{\mathbbm{H}^N}}}{\text{osc}}u(\xi)},
\end{align}
which immediately offers
\begin{align}
	\displaystyle{\underset{{B_{2^{-k-1}}^{\mathbbm{H}^N}}}{\text{osc}}u(\xi)}&\leq {(1-\delta)}\displaystyle{\underset{{B_{2^{-k}}^{\mathbbm{H}^N}}}{\text{osc}}u(\xi)}\\
	&\leq 2(1-\delta)^{k+1}\\
	&=2.2^{-\c(k+1)}.
\end{align}
Let for some $k>0,\, \xi\in B^{\mathbbm{H}^N}_{2^{-k}}\setminus B^{\mathbbm{H}^N}_{2^{-(k+1)}}.$ This grants
\begin{align}
	|u(\xi)-u(0)|\leq 2.2^{-\c k}=2.2^\c2^{-\c(k+1)}\leq C|\xi|_{\mathbbm{H}^N}^\c,
\end{align}
for $C=2.2^\c.$ Hence the claim. \qed

\section{Funding and/or Conflicts of interests/Competing interests}
Both authors thank DST/SERB for the financial support under the grant CRG/2020/000041.\\

There are no conflict of interests of any type. This manuscript does not use any kind of data.

\end{document}